\documentclass{article}

\usepackage{graphicx,url,xcolor}
\usepackage{caption,subfigure}
\usepackage{amsfonts,amssymb,amsmath,amsthm}
\newtheorem{example}{Example}
\newtheorem{definition}{Definition}
\newtheorem{proposition}{Proposition}
\newtheorem{appli}{Application}

\def\ind{{\mathchoice {\rm 1\mskip-4mu l} {\rm 1\mskip-4mu l}{\rm 1\mskip-4.5mu l} {\rm 1\mskip-5mu l}}}

\begin{document}

 \author{Boualem Djehiche, Alain Tcheukam and  Hamidou Tembine
\thanks{B. Djehiche is with Department of Mathematics, KTH Royal Institute of Technology, Stockholm, Sweden }\thanks{
 A. Tcheukam and H. Tembine are with Learning and Game Theory Lab, New York University Abu Dhabi. Email: tembine@nyu.edu}
 }
  %\ead{tembine@nyu.edu}

\title{Mean-Field-Type Games in Engineering}


%\author{}

%\address{}
\maketitle
\begin{abstract}
A mean-field-type game is a  game in which the instantaneous payoffs and/or the state dynamics  functions involve  not only the state and the action profile but also the joint distributions of state-action pairs.This article presents engineering applications of mean-field-type games. 
\end{abstract}


\maketitle

\section{Introduction}
With the ever increasing amounts of data becoming available, strategic data analysis and decision-making will become  more pervasive as a necessary ingredient for societal infrastructures.  
In many network engineering games, the performance metrics  depend on some few aggregates of the parameters/choices. 
A typical example  is the congestion field in traffic engineering where classical cars and smart autonomous driverless cars  create   traffic congestion levels on  the roads. 
The congestion field can be learned, for example by means of crowdsensing,  and can be used for efficient and accurate prediction of the end-to-end delays of commuters.
Another  example is the interference field where it is the aggregate-received signal of the other users that matters rather than their individual input signal.  In such games, in order for a transmitter-receiver pair to determine his best-replies, it is unnecessary that the pair is informed about the other users' strategies. If a user is informed about the aggregative terms given  her own strategy, she will be able to efficiently exploit such  information to perform better.  In these situations the outcome is influenced not only by the state-action profile but also by the distribution of it. The interaction can be captured by a game with distribution-dependent payoffs called {\it mean-field-type games} (MFTG).  An MFTG is   basically a game in which the instantaneous payoffs and/or the state dynamics  functions involve  not only the state and the action profile of the agents but also the joint distributions of state-action pairs.

{\color{black} The main contributions of this article can be summarized as follows. The first contribution of this article is the review of some relevant engineering applications of MFTG.
  Considering Liouville type systems with drift, diffusion and jumps, that are dependent on  time-delays,
 state mean-field and action mean-field terms. Proposition \ref{mainr1} establishes an equilibrium equation for non-convex action spaces. Proposition \ref{mainr2}  provides  a  stochastic maximum principle   that covers  decentralized information  and partial observation systems which are crucial in engineering systems. Various engineering applications in discrete or continuous  variables (state, action or time) are provided. Explicit solutions are provided in propositions \ref{onlinemeeting:res1}  and \ref{delay:res2} which are mean-field type game problems with non-quadratic costs.

  The article is structured as follows. The next section overviews  earlier works on  static mean-field games, followed by discrete time mean-field games with measure-dependent  transition kernels.  Then, a basic MFTG with finite number of agents is presented. After that, the discussion is divided into two illustrations in each of the following  areas of engineering (Figure \ref{fig:content}) :   Civil Engineering (CE), Electrical Engineering (EE), Computer Engineering (CompE), Mechanical Engineering (ME), General Engineering (GE). 
\begin{itemize}\vspace{-3mm}\item CE:  road traffic networks  with random incident states and multi-level building evacuation
\item EE:   Millimeter wave wireless communications and  distributed power networks
 \item CompE:  Virus spread over networks and  virtual machine resource management in cloud networks
   \item  ME:   Synchronization of oscillators, consensus, alignment and  energy-efficient buildings
   \item GE:    Online meeting: strategic arrivals and starting time and   mobile crowdsensing as a public good.
   \end{itemize}
\vspace{-3mm}   
The article proceeds by presenting the effect of time delays of coupled mean-field dynamical systems and decentralized information structure. Then, a discussion on the drawbacks, limitations, and challenges of MFTGs is highlighted. Lastly, a summary of the article and concluding remarks are presented. 
\begin{figure}
\centering
\includegraphics[width=0.4\textwidth]{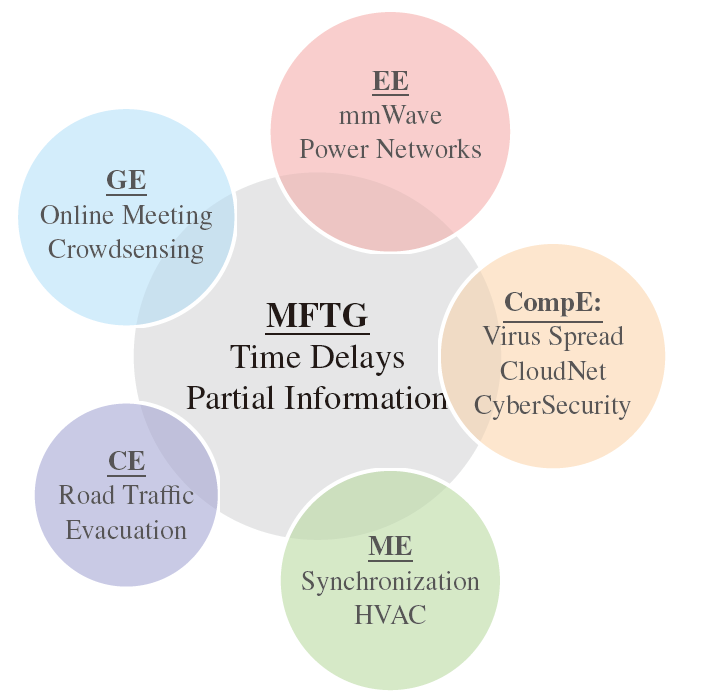}
  \caption{MFTG with engineering applications covered in this work.} \label{fig:content}
  \end{figure}

\subsection{Mean-Field Games: Static Setup}
A static mean-field game is one in which all users make choices (or select a strategy) simultaneously, without knowledge of the strategies that are being chosen by other users and the game is played once.
Any mean-field game with sequential moves is a dynamic mean-field game. In this work, games which are played more than once will be considered as dynamic game.  
This subsection overviews static mean-field games and games in which the underlying processes are in stationary regime (time-independent).
 Mean-field games have been around for quite some time in one form or another, especially in transportation networks and in competitive economy. 
In the context of competitive market with large number of agents, a 1936 article \cite{wald} captures the assumption made in mean-field games with large number of agents, in which the author  states: 

``each of the participants has the opinion that its own actions do not influence the prevailing price".

Another comment on the impact on the population mean-field term was given in \cite{VM1} page 13: 

`` When the  number of participants becomes large, some hope emerges  that of the influence of every particular participant will become negligible \ldots"

The population interaction involves many  agents for each type or class and location, a common approach is to replace the
individual agents' variables and to use continuous variables to represent the aggregate average of type-location-actions. In the large population regime, the mean field limit is then modeled by state-action and location-dependent time process (see Figure \ref{fig:mfg}). This type of aggregate models are also known as non-atomic or population games. It is closely related to the mass-action interpretation in \cite{nash51},  Equation (4) in page 287.

%%\end{document}

\begin{figure}
\centering
    \includegraphics[height=4cm,width=6cm]{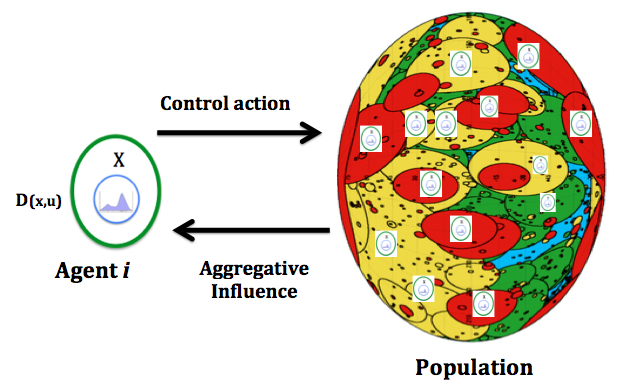} 
\caption {Each agent  is with its  own state and own mean-field interacts with the aggregates from  the rest of the population. The population mean-field  is formed from the reaction of the agents and affects the behavior of the individual agents and their own mean-field. } \label{fig:mfg}
\end{figure}

  In the context of transportation networks, the mean-field game 
 framework, underlying the key foundation,  goes back to the pioneering works of \cite{wardrop} in the 1950s.   
Therein, the basic idea is to describe and understand interacting traffic flows among a large population of  agents moving from multiple sources to destinations, and interacting with each other. The congestion created on the road and at the intersection are subject to capacity and flow constraints.  This corresponds to a constrained mean-field game problem as noted in \cite{beckmann}. 
A common behavioral assumption in the study of transportation and communication networks is that travelers or packets, respectively, choose routes that they perceive as being the shortest under the prevailing traffic conditions.  As noted in  \cite{Knight},    collection of individual decisions may result to a situation which drivers cannot reduce their journey times by unilaterally choosing another route. 
The work in \cite{Knight} such a  resulting traffic pattern as an equilibrium. Nowadays, it is indeed known as the Wardrop  equilibrium \cite{wardrop,Dafermos}, and it is  thought 
of as a steady state obtained after a transient phase in which travelers successively adjust their route choices 
until a situation with stable route travel costs and route flows has been reached \cite{congest2,conges}. In the seminal contribution \cite{wardrop}, page 345 the author stated two principles that formalize this notion of 
equilibrium and the alternative postulates of the minimization of the total travel costs. His first principle reads:

{``The journey times on all the routes actually used are equal, and less than those which would be experienced by a single vehicle on any unused route."}

Wardrop's first principle of route choice, which is identical to the notion postulated in \cite{Kohl,Knight}, became widely used as a sound and simple behavioral principle to describe the 
spreading of trips over alternate routes due to congested conditions. Since its introduction in the context of transportation networks in 1952 and its mathematical formalization by \cite{beckmann,smith}  transportation planners have been using Wardrop equilibrium models to predict commuters decisions in real-life networks.  

{\color{black}
The key congestion factor is the flow or the fraction of travelers per edge on the roads (see Application 1). The above Wardrop problem is indeed a mean-field on a discrete space. The exact mean-field term here corresponds to a mean-field of actions (a choice of a route). Putting this in the context of infinite number of commuters results to  end-to-end travel times that are function of own choice of a route and the mean-field distribution of travelers across the graph (network). 
}

In a population context, the equilibrium concept of  \cite{wardrop}  corresponds to a 
Nash equilibrium of the mean-field game with infinite number of agents. 
The works \cite{Dafermos,Dafermos2} provide a variational formulation of  the (static) mean-field equilibrium.  
 
 The game theoretic models such as evolutionary games \cite{bill1,bill2,bill3}, global games \cite{global1,global2}, anonymous games, aggregative games \cite{dubey}, population games \cite{bill4,bill5,bill6}, and large games, share several common features. Static mean-field games with large number of agents were widely investigated (see \cite{conti2,conti5,conti1,conti3,conti4,aumann} and the references therein).

\subsection{Mean-Field Games: Dynamic Setup}
This section overviews mean-field games which are dynamic (time-varying and played more than once) and their applications in engineering.
{\color{black}  
\vspace{-2mm}
\begin{definition}[Mean-Field Game: Infinite Regime] \label{defi01}
\emph{A (homogeneous population) mean-field game (MFG) is a  game in which the instantaneous payoff of a generic agent (say 1) and/or the state dynamics coefficient functions involve an individual  state-action pair $x_{1t},u_{1t}$  and the distribution of state-action pairs of the other decision-makers, $m_t$ at time $t.$  The individual state and action spaces are identical across the homogeneous population denoted by $\mathcal{X}_i=\mathcal{X}_1, U_j=U_1$ for all $i.$ The state transition to the next state follows $ \mathbb{P}( . | x_{1t}, u_{1t}, m_t).$
Thus, the instant payoff  function of   a generic agent  (say $j$) has the following structure:} $$r_i = r_1 = r:\ \mathcal{X}_1\times U_1\times \mathbb{P}( \mathcal{X}_1\times U_1) \rightarrow \mathbb{R},$$ \emph{with} $r(x_{1t}, u_{1t}, m_t).$ 
\end{definition}
\vspace{-3mm}
The mean-field game model has been extended to include several other features such as incomplete information, common noise, heterogeneous population, finite population or a mixture between  finite number of clusters and infinite population regimes. 

}

The key ingredients of  dynamic mean-field games appeared in  \cite{jova,jova2} in the early 1980s.  
The  work in  \cite{jova} proposes a game-theoretic model that explains why smaller firms grow faster and are more likely to fail than larger firms in large economies. 
The game is played over a discrete time space. 
Therein, the mean-field is the aggregate demand/supply which generates a price dynamics. The price moves forwardly, and the agents react to the price and generate a demand and the firm produces a supply 
with associated cost, which regenerates the next price and so on. The author introduced a  backward-forward system to find equilibria (see for example Section 4, equations D.1 and D.2 in \cite{jova}). 
The backward  equation is obtained as an optimality to the individual response, i.e., the value function associated with the best response to price, and the forward equation for the evolution of price. Therein, the consistency check is about the mean-field of equilibrium actions (population or mass of actions), that is, the equilibrium price solves a fixed-point system: 
the price regenerated after the reaction of the agents through their individual best-responses should be consistent with the price they responded to. 

Following that analogy, a more general framework was developed in \cite{jova2}, where the mean-field equilibrium is introduced in the context of dynamic games with large number of decision-makers. A mean-field equilibrium is defined in  \cite{jova2}, page 80  by two conditions: (1) each generic agent's action is best-response to the mean-field, and (2) the mean-field is consistent and is exactly reproduced from the reactions of the agents. This matching argument was widely used in the literature as it can be interpreted as a generic agent reacting to an evolving mean-field object and at the same time the mean-field is formed from the contributions of all the agents. The authors of \cite{bergins} show  how common noise can be introduced into the mean-field game model  (the mean-field distribution evolves stochastically) and extend the {\color{black} Jovanovic-Rosenthal existence theorem \cite{jova2}.}

Continuous time version of  the works  \cite{jova,jova2} can be found 
in \cite{benamou,benamou2,PLref3,peter}. The reader is referred  to \cite{alainv1,carmona,gomes,wil,gomes2t,gomes3t,gomes4t}  for recent development of mean-field game theory. 
 The authors \cite{PLref1,PLref2,PLref3,PLref4,PLref5,PLref6} have developed  a powerful tool for modelling strategic behavior of large population  of agents, each of them having a negligible impact on the population mean-field term. 
Weak solutions of mean-field games are analyzed in \cite{porreta}, Markov jumps processes \cite{wang,wang2}, and leader-followers models in \cite{team4}. 
Finite state mean-field game models were analyzed in \cite{weintraub1,weintraub2,tembinegamenets,finitediogo2,finitediogo1,finitediogo1b,finitediogo1c,johari}.  Team and social optimum solutions can be found in \cite{team1, team2,team3,team4,k11}. The  work in \cite{sznitman,moon3,mini} provide mean-field convergence  of  a class of McKean-Vlasov dynamics.  Numerical methods for mean-field games can be found in \cite{yves,yves2,yves3,yves4}.

%\subsection*{ Relevance of MFTGs for engineering applications: }
{\color{black} Table \ref{table:ref} summarizes some engineering applications of mean-field-type game models. }

 \begin{table}
 \centering
\vspace{-24.55mm}
\caption{Some applications of  MFTGs in Engineering.  HVAC stands for (heating, ventilation and air conditioning) systems. D2D stands for Device-to-Device (D2D) }
\begin{tabular}{ll}
\hline
Area & Works\\ \hline
planning & \cite{planning} \\  
state estimation and filtering &  \cite{filter,sergio}\\  
synchronization & \cite{pras1,pras2,pras3,pras4}\\  
opinion formation & \cite{pe9}\\ 
network security& \cite{k10,sec,moon,moon2,ccdc16two} \\  
power control & \cite{powermi,wireless1,wireless5}\\  
medium access control & \cite{wireless2,evo1}\\  
cognitive radio networks &  \cite{wireless3, wireless4}\\  
electrical vehicles & \cite{romain,fe8}\\  
scheduling & \cite{sha1}\\ 
cloud networks & \cite{sha2,cloudnet1,cloudnet2}\\  
wireless networks & \cite{wir}\\ 
auction & \cite{auction1,auction2}\\ 
cyber-physical systems & \cite{cyber,small3}\\ 
airline networks &  \cite{ssdv1}\\  
sensor networks & \cite{aziz}\\ 
traffic networks &\cite{filter1mf,filter2mf,nyue,pri7} \\ 
big data & \cite{bdata}\\ 
%small cell networks &  \cite{small1,small4}\\ 
D2D networks & \cite{small2,psy,pri8}\\ 
multilevel building evacuation & \cite{aims,burj2,burj3,burj4}\\ 
power networks & \cite{pe1,pe2,pe3,pe4,pe5,pri4,pri6} \\ 
& \cite{pe6,pe7,pe8,fe8,fe9} \\ 
& \cite{fe10,ccdc16one,ye}\\  
HVAC & \cite{pe10, pe11,pe12,pe13,pe14,tac1}\\ \hline
\end{tabular}
\label{table:ref}
\end{table}

 \subsubsection{ Limitations of the existing mean-field game models}

Most of the  existing mean-field game models share the following assumptions:

 Big size: A typical assumption is to consider   an infinite number decision-makers, sometimes,  a continuum of decision-makers. 
The idea of a continuum of decision-makers may seem outlandish to the reader. Actually, it is no stranger than a continuum  of particles used in fluid mechanics, in water distribution, or in petroleum engineering.
In terms of practice  and experiment however, decision-making problems with continuum of decision-makers is rarely observed in engineering.
There  is a huge  difference between a fluid with a continuum of particles and a decision-making problem with a continuum of agents. Agents may  physically occupy a space  (think of agents inside a building or a stadium) or a resource, and the size or number of agents that most of 
engineering systems can handle can be relatively large or growing but remain currently finite \cite{smc}. It is in part due to  the limited resource per shot or limited number of servers at a time. In all the examples and applications  provided below, we still have a finite number of interacting agents.  Thus, this assumption appears to be very restrictive in terms of engineering applications.

Anonymity: The index of the decision-maker does not affect the utility. The agents are assumed to be indistinguishable within the same class or type.
The drawback of this assumption is that most individual  decision-makers in engineering are in fact not  necessarily anonymous 
(think of Google, Microsoft, Twitter, Facebook, Tesla), the classical mean-field game model is inappropriate, and does not apply to such  situations.  In mean-field games with several 
types (or multi-population mean-field games), it is still assumed  that there is large number of agents per type/class/population, which is not realistic in most of the engineering applications considered in this work.

 NonAtomicity: A single decision-maker has a negligible effect on the mean-field-term and on the global utility. 
One typical example where this assumption  is not satisfied is a situation of targeting a room comfort temperature, in which the air conditioning controller adjusts the heating/cooling depending on the temperature in the room, the temperatures of the other connecting zones and the ambient temperature. It is clear  that the decision of the controller to heat or to cool affect the variance of the temperature inside the room. Thus, the effect of the individual action of that controller 
on the temperature distribution (mean-field) inside the room cannot be neglected.

To summarize, the above conditions appear to be very restrictive in terms of engineering applications, and to overcome this issue
a more flexible  MFTG framework has been proposed.

\subsubsection{What MFTGs can bring to the existing decision-making models?}

MFTGs not only relax of the  above assumptions but also incorporate the behavior of the agents as well as their effects in the  mean-field terms and in the outcomes (see Table  \ref{table:reff1}). 
{\color{black} 
\begin{itemize}
\vspace{-4mm}
\item[(1)]
In MFTGs, the number of users can be finite or infinite.
\item[(2)]
The indistinguishability property (invariance in law by permutation of index of the users) is not assumed in MFTGs. 
\item[(3)] A single user may have a non-negligible impact of the mean-field terms, specially in the distribution of own-states and own mixed strategies.
\end{itemize}
\vspace{-4mm}
These  properties (1)-(3) make  strong differences between mean-field games and MFTGs (see \cite{bou3} and the references therein). 
}

 \begin{table}
 \centering
\caption{Key limitations and differences between the game models}
\begin{tabular}{llll}
\hline
Area & Anonymity & Infinity &  Atom\\ \hline
population games \cite{wardrop,beckmann} & yes &  yes & no\\  
evolutionary games \cite{maynard} &  yes &  yes & no\\  
non-atomic games \cite{jova2} & yes &  yes & no\\ 
aggregative games \cite{dubey}&  &  relaxed &  \\  
global games \cite{global1,global2}& yes &  yes & no\\  
large games \cite{conti2}& yes &  yes & no\\  
anonymous games \cite{jova2}& yes &  yes & no\\  
 mean-field games&  yes &  yes & no\\  
 nonasymptotic mean-field games&  nearly&  no & yes\\  
MFTG & relaxed & relaxed & relaxed\\  \hline
\end{tabular}
\label{table:reff1}
\end{table}

MFTG seems to be more appropriate in such engineering situations because it does not assume indistinguishability, it captures the effect of each agent
 in the distribution and the number of agents is arbitrary 
as we will see below.   Table \ref{table:notations} summarizes the notations used in the manuscript.

\begin{table}[htbp]\caption{Table of Notations}
\centering % to have the caption near the table
\begin{tabular}{r c p{10cm} }
%\toprule 
\hline
$\mathcal{I}$ &  $\triangleq$ &set of decision-makers\\
$T$ & $\triangleq$ & Length of the horizon\\
$[0,T]$ & $\triangleq$ & horizon of the mean-field-type game\\
$t$ & $\triangleq$ & time index\\
$\mathcal{X}$ & $\triangleq$ & state space\\
$W$ & $\triangleq$ & Brownian motion\\
$\sigma$ & $\triangleq$ & Diffusion coefiicient\\
$N$ & $\triangleq$ & Poisson jump process\\
$\gamma$ & $\triangleq$ & Jump rate coefiicient\\
${U}_i$ & $\triangleq$ & control action  space of agent $i\in \mathcal{I}$\\
$\mathcal{U}_i$ & $\triangleq$ & admissible strategy  space\\
$u_i$ & $\triangleq$ & state space\\
$r_i$ & $\triangleq$ & instantaneous payoff\\
$D_{(x,u)}$& $\triangleq$ & distribution of state-action\\
$R_i$ & $\triangleq$ & Long-term payoff functional\\
\hline
\end{tabular}
\label{table:notations}
\end{table}

\subsection{Background on MFTGs}
This  section presents a background on  MFTGs.  
\vspace{-2mm}
\begin{definition}[Mean-Field-Type Game] \label{defi1}
\emph{A mean-field-type game (MFTG) is a  game in which the instantaneous payoffs and/or the state dynamics coefficient functions involve  not only the state and the action profile but also the joint distributions of state-action pairs (or its marginal distributions, i.e., the distributions of states  or the distribution of actions). Let $\mathcal{I}$ be the set of agents, $\mathcal{X}_i$   the state  space of agent $i$   and  $ \mathcal{X}:=\prod_{i\in \mathcal{I}}\mathcal{X}_i=\mathcal{X}_1\times \mathcal{X}_2\times \ldots $  the state profiles space of all  agents. $U_i$ is the action space of agent $i$ and $U=\prod_j U_j $ is the action profile space of all agents.
A typical example of payoff  function of agent  $j$ has the following structure: $$r_i:\ \mathcal{X}\times U\times \mathbb{P}( \mathcal{X}\times U) \rightarrow \mathbb{R},$$ with $r_i(x,u, D_{(x,u)})$ where $(x,u)$ is the state-action profile of the agents and $D_{(x,u)}$ is the distribution of the state-action pair $(x,u).$  $\mathcal{X}$ is the state space,  $U$ is the action profile space of all agents and $\mathbb{P}( \mathcal{X}\times U)$ is the set of probability measures over $\mathcal{X}\times U .$ }
\end{definition}

\vspace{-2mm}
From Definition \ref{defi1}, a mean-field-type  game can be static or dynamic in time.
One may think that MFTG is  a small and particular class of games. However, this class includes the classical games in strategic form because any payoff function $r_i(x, u)$ can be written as $r_i(x,u, D).$ 
 
 When randomized/mixed strategies are used in the von Neumann-type payoff, the resulting payoff can be written as $E[r_i(x,u)]=\int r_i(x,u) D_{(x,u)}(dx,du)=\hat{r}_i(D).$ Thus, the form $r_i(x,u, D)$ is more general and includes non-von Neumann payoff functions.
\vspace{-2mm} 
  \begin{example}[Mean-variance payoff] \emph{The payoff function of agent $i$ is $E[r_i(x,u)]-\lambda \sqrt{var[r_i(x,u)]}, \lambda\in \mathbb{R}$ which can be written as a function of $r_i(x,u,D_{(x,u)}).$ {\color{black} For any number of interacting agents, the term $D_{(x_i,u_i)}$ plays a non-negligible role in the standard deviation $\sqrt{var[r_i(x,u)]}.$ Therefore, the impact of agent $i$ in the individual mean-field term   $D_{(x_i,u_i)}$ cannot be neglected.}}
\end{example}

\vspace{-2mm}
\begin{example}[Aggregative games] 
\emph{The payoff function of each agent depends on its own action and an aggregative term of the other actions. Example of payoff functions include  $r_i(u_i, \sum_{j\neq i} u^{\alpha}_j), \ \alpha>0$ and 
 $r_i(x_iu_i, \sum_{j\neq i} x_ju_j).$ }
\end{example}

\vspace{-2mm}
 In the non-atomic setting, the influence of an individual state $x_i$ and individual action $u_i$ of any agent $i$ will have a negligible impact on mean-field term  $\hat{D}_{(x,u)}=\lim_{n\rightarrow +\infty}\ \frac{1}{n-1}\sum_{j\neq i}\delta_{\{ x_j,u_j\}}.$ In  that case, one gets to the so-called mean-field game. 
\vspace{-2mm}
\begin{example}[Population games]\emph{ Consider a large population of agents. Each agent has a certain state/type $x \in \mathcal{X}$ and can choose a control action $u\in \mathcal{U}(x).$ Let the proportion of type-action of the population as $m.$
The payoff of the agent with type/state $x,$ control action $u$ when the population profile $m$ is $r(x, u,m).$  
{\it Global games} with continuum of agents were studied in \cite{global1} based on the Bayesian games of \cite{global2}, which uses the proportion of actions.}
\end{example}
\vspace{-2mm}
In the case where  both non-atomic and atomic terms are involved in the payoff, one can write the payoff  as $r_i(x,u, D, \hat{D})$ where $\hat{D}$ is the population state-action measure. 
 Agent $i$ may  influence $D_i$ (distribution of its own state-action pairs) but  its influence on $ \hat{D}$  may be limited.

The main goals of static mean-field-type games are: (1) identify solution concepts \cite{solant} such as Nash equilibrium, Bayesian equilibrium, correlated equilibrium, Stackelberg solution etc. (2) Computation of solution concepts. (3) Development of algorithms and learning procedures to reach and select efficient equilibria, (4) Mechanism design for incentivizing agents.
 The next section presents a class of dynamic MFTGs which are played over several stages.

\section{A Basic Dynamic MFTG: Finite Regime}
Consider a basic MFTG  with $n\geq 2$ agents interacting over horizon $[0,T],\ \ T>0. $ The individual state dynamics of agent $i$ is given by
\begin{equation}\label{eq:state} 
 \begin{array}{ll}
 dx_i=b_i\left(x_i,u_i, D_{(x_i,u_i)}, \frac{\sum_{k\neq i}\delta_{(x_k,u_k)}}{n-1}\right) dt+ \sigma_i\left(x_i,u_i, D_{(x_i,u_i)}, \frac{\sum_{k\neq i}\delta_{(x_k,u_k)}}{n-1}\right)  dW_i,\\  x_i(0) \sim D_{i,0}\end{array}
 \end{equation}
and the payoff functional of agent $i$ is 
\begin{eqnarray}\label{eq:reward} R_i(u)=g_i\left(x_i(T), D_{x_i(T)}, \frac{\sum_{k\neq i}\delta_{x_k(T)}}{n-1}\right)  +\int_0^T r_i\left(x_i,u_i, D_{(x_i,u_i)}, \frac{\sum_{k\neq i}\delta_{(x_k,u_k)}}{n-1}\right)  dt,
\end{eqnarray}
where  the strategy profile is $u=(u_1,\ldots, u_n),$ which also denoted as $(u_i,u_{-i}).$ The functions $b_i, \sigma_i,g_i, r_i$ are  measurable functions.  $x_i(t): = x_i(t)[u]$ is the state of agents $i$ under of the strategy profile $u,$ $ D_{x_i(t)}=\mathcal{L}(x_i(t))$ is the probability distribution (law) of $x_i(t).$  $ D_{(x_i(t),u_i(t))}=\mathcal{L}(x_i(t),u_i(t))$ is the probability distribution of the state-control action pair $(x_i(t),u_i(t))$ of agent $i$ at time $t.$  $\delta_y$ is the $\delta-$Dirac measure concentrated at $y,$  and $W_i$ is a standard Brownian motion defined over the filtration  $(\Omega, \mathbb{P}, (\mathcal{F}_t)_{t\leq T}).$ 

The novelty in the modelling  of (\ref{eq:state})-(\ref{eq:reward})  is that  each individual agent $i$ 
influences  its own mean-field terms $D_{x_i(t)},$ and  $ D_{(x_i(t),u_i(t))}$ independently on the total number of interacting agents. In particular, the influence of agent $i$ on those mean-field terms remain non-negligible even when there is a continuum of agents.
The distributions $D_{x_i}$ and $D_{(x_i,u_i)}$ represent two important  terms in the modeling of MFTGs. These terms are referred to as {\it individual mean-field} terms. In the finite regime, the other agents are captured by the empirical measures $\frac{\sum_{k\neq i}\delta_{x_k}}{n-1}$ and $ \frac{\sum_{k\neq i}\delta_{(x_k,u_k)}}{n-1}.$  We refer  these terms to as {\it population mean-field} terms. 

Similarly, a basic  discrete time (discrete or continuous state) MFTG  is given by 
\begin{equation}\label{eq:discrete}  \begin{array}{ll}
 x_{i,t+1}\sim q_i\left(.| x_{i,t},u_{i,t}, D_{(x_{i,t},u_{i,t})}, \frac{\sum_{k\neq i}\delta_{(x_{k,t},u_{k,t})}}{n-1}\right) ,\\  x_{i0} \sim D_{i,0} \\
 R_i(u)=g_i\left(x_{i,T}, D_{x_{i,T}}, \frac{\sum_{k\neq i}\delta_{x_{k,T}}}{n-1}\right) +\sum_{t=0}^{T-1} r_i\left(x_{i,t},u_{i,t}, D_{(x_{i,t},u_{i,t} )}, \frac{\sum_{k\neq i}\delta_{(x_{k,t},u_{k,t})}}{n-1}\right),
\end{array}
\end{equation}
where $q_i(.|.)$ is the transition kernel of agent $i$ to next states. 

Mean-field-type control and global optimization can be found in \cite{andersson,alainv1,forna,pri2,pri3,pri5}. The models (\ref{eq:state}) and  (\ref{eq:discrete}) are easily adapted to bargaining solution, cooperative and coalitional MFTGs and can be found  in \cite{coalition,duncan,duncan2}. Psychological MFTG was recently introduced in \cite{psy,duncan3} where spitefulness, altruism, selfishness, reciprocity of the agents are examined by means empathy, other-regarding behavior and psychological factors.

\vspace{-2mm}
\begin{definition}
\emph{An admissible control strategy of agent $i$ is an $\mathcal{F}_i-$adapted and square integrable process with values in a non-empty subset ${U}_i.$   Denote by $\mathcal{U}_i =L^2_{\mathcal{F}_i}([0,T],\ U_i)$  the class of admissible control strategies of agent $i$.} 
\end{definition}

\vspace{-3mm}
\begin{definition}[Best response]
\emph{Given a strategy profile of the other agents $(u_1,\ldots, u_{i-1},u_{i+1},\ldots, u_n),$  with $u_j,\ j\neq i$ that are admissible and the mean-field terms $D$,  the best response problem of agent $i$ is: }
\begin{eqnarray}
\label{smp0best}
\left\{
\begin{array}{l}
\displaystyle{\sup_{u_i\in \mathcal{U}_i} \mathbb{E} \left[ R_i(u)\right] },\\
\mbox{subject  to} \  (\ref{eq:state}) 
\end{array} \right.
\end{eqnarray}
\end{definition} 

%has at least one solution in the space of square integrable function over $[0,T]\times D \rightarrow \mathbb{R}^m.$ 

The first goal is to find and characterize the best response strategies of each user. For user $i$ it consists to solve problem (\ref{smp0best}).
In problem (\ref{smp0best}), the information structure that is available to user plays in important role. We will distinguish three type of strategies: (1) open-loop strategies that are only measurable function of $t,$ (2) state-feedback strategies that are measurable functions of state and time, (3) state-and-mean-field feedback strategies that measurable functions of state, mean-field and time.  To solve problem (\ref{smp0best}), four different methods have been developed:
\begin{itemize}
\vspace{-3mm}
\item Direct approach which consists to write the payoff functional in a form such that the optimal value and optimizers are trivially obtained, and a verification and validation procedure follows.
\item A stochastic maximum principle (Pontryagin's approach) which provides necessary conditions for optimality. 
\item A dynamic programming principle (Bellman's approach) which consists to write  the value of the problem (per agent) in  (backward) recursion form, or as solution to a dynamical system.
\item Uncertainty quantification approach by means of Wiener chaos expansion of all the stochastic terms and the use of Kosambi-Karhunen-Loeve expansion which is a representation of a stochastic process as an infinite linear combination of orthogonal functions, analogous to a Fourier series representation of a function over a bounded domain.
\end{itemize}

\vspace{-3mm}
If every user solves its best-response problem, the resulting system will be a Nash equilibrium system defined below.
\vspace{-2mm}
\begin{definition} 
\emph{A (Nash) equilibrium of the  game is a strategy profile $(u^*_1,\ldots, u^*_n)$ such that  for every agent $i$,  $$\mathbb{E}[R_i(u^*)]\geq \mathbb{E}[R_i(u^*_1,\ldots, u_{i-1}^*,u_{i},u_{i+1}^*,\ldots, u^*_n)], $$ for all   $ u_i\in \mathcal{U}_i .$}
\end{definition} 
 The second goal is to find and characterize Nash equilibria of the mean-field-type game.  We provide below a basic example in which  the Nash equilibrium problem can be solved semi-explicitly using Riccati system.
\vspace{-2mm}
\begin{example}[Network Security Investment \cite{k10} ] \label{example:security}
\emph{A graph is connected if there is a path that  joins any point to any other point in the graph. Consider  $n\geq 2$  decision-makers over a connected graph.
Thus, the security of a node is influenced by the others through possibly multiple hops. The effort of user $i$  in security investment  is $u_i.$ The associated cost  may include  money (e.g., for purchasing antivirus software), time and energy (e.g., for system scanning, patching).  Let $x(t)$ be the security level of the network at time $t$ and}
\begin{equation} \begin{array}{lll} 
R_i(u)= \  \ - \frac{1}{2}[x(T)-Ex(T)]^2  
+\int_0^T q_i(t)x(t)(1-\epsilon_i(t) x(t))-\rho_i(t) u_i(t)-\frac{r_i(t)}{2}u^2_i(t)dt.
\end{array}\end{equation}  \emph{The best-response  of  user  $i$ to $(u_{-i}, E[x]): = (u_1,\ldots, u_{i-1},u_{i+1},\ldots, u_n, E[x]),$ solves the following   linear-quadratic mean-field-type control problem }\begin{eqnarray}
\label{LQ00game1type} 
\left\{
\begin{array}{lll} \sup_{u_i\in \mathcal{U}_i}  E \left[ R_i(u_1,\ldots, u_n) \right],
\
\displaystyle{\mbox{ subject to }\ }\\
dx=\left\{-ax-\bar{a}E[x]+\sum_{i=1}^n b_i u_i\right\}dt + cxdW, \\
x(0)\in \mathbb{R}, 
 \end{array}
\right.
\end{eqnarray}
\emph{where, $q_i(t)\geq 0, \epsilon_i(t)\geq 0, \rho_i(t)\geq 0,\  r_i(t)>0$ and  $a,\bar{a}, b_i, c$ are real numbers  and where $E[x(t)]$ is the 
expected value of network security level created by all users  under the control action profile $(u_1,\ldots, u_n).$ Note that the expected 
value of the terminal term in $R_i$  can be seen as a weighted  variance  of the state \cite{tac1} since 
$E[(x(t)-E[x(t)])^2]= var(x(t)).$ The optimal control action   is in state-and-mean-field feedback  form:}
\begin{equation} \label{LQ-ucontroltype}\begin{array}{ll} u^*_i(t)=-\frac{b_i}{r_i(t)}\left[ \beta_i(t) x(t)+\eta_{1i}(t) E[x(t)]+\eta_{2i}(t)\right]-\frac{\rho_i(t)}{r_i(t)}, \nonumber \\
0=\dot{\beta}_i +(-2a+c^2)\beta_i -\beta_i\sum_{j=1}^n\frac{b^2_j}{r_j}\beta_j +2q_i\epsilon_i,\\ 
\beta_i(T)=1,\\
\dot{\eta}_{1i}-2(a+\bar{a})\eta_{1i} -2\bar{a}\beta_i- \beta_i \sum_{j=1}^n\frac{b^2_j}{r_j}\eta_{1j} -\eta_{1i} \sum_{j=1}^n\frac{b^2_j}{r_j}(\beta_j+\eta_{1j}) =0,\\ 
 {\eta}_{1i}(T)=-1,\\
\dot{\eta}_{2i}- (a+\bar{a})\eta_{2i} - \beta_i \sum_{j=1}^n\frac{b_j}{r_j}  (b_j\eta_{2j}+\rho_j)  -\eta_{1i} \sum_{j=1}^n \frac{b_j}{r_j}(b_j\eta_{2j}+\rho_j) -q_i =0,\\
{\eta}_{2i}(T)=0.
\end{array}
\end{equation}

\emph{Figure \ref{rnstateneu} plots  the optimal cost trajectory with the step size  $2^{-8},$ the horizon is $[0,1],$ the other parameters  are $b=5, r=1, q=1, \rho=0.0001, \epsilon=0.1.$
Figure \ref{rnstateneu3} plots the optimal state vs the equilibrium state.  As noted in \cite{xinluo},  the security state is higher when there is a cooperation between the users and when the coalition formation cost is small enough.  The inefficiency of Nash equilibria behavior is widely known in game theory in which the Nash equilibrium can be inefficient compared to the global optimum of the system. The relative payoff difference between the worse Nash equilibrium payoff and global optimum payoff have been proposed in the literature \cite{poa1,tanimoto} as measure of inefficiency. Another measure of inefficiency is the Price of anarchy, which has been proposed in \cite{christos1,christos2}. It  measures the ratio between the worse Nash equilibrium payoff and global optimum payoff. Note however that, one needs to be careful by taking a ratio here, because the denominator may vanish in our context. Note that restricting the analysis to the set of symmetric strategies may lead to performance degradation \cite{pri3} as symmetric Nash equilibria may not be performant even in symmetric games. Thus, looking at $\epsilon-$Nash equilibria via mean-field limiting behavior does not help in improving the efficiency of  Nash equilibria. }

\begin{figure}[htb]
\centering
\begin{minipage}{0.9\linewidth}
 \includegraphics[width=0.8\textwidth]{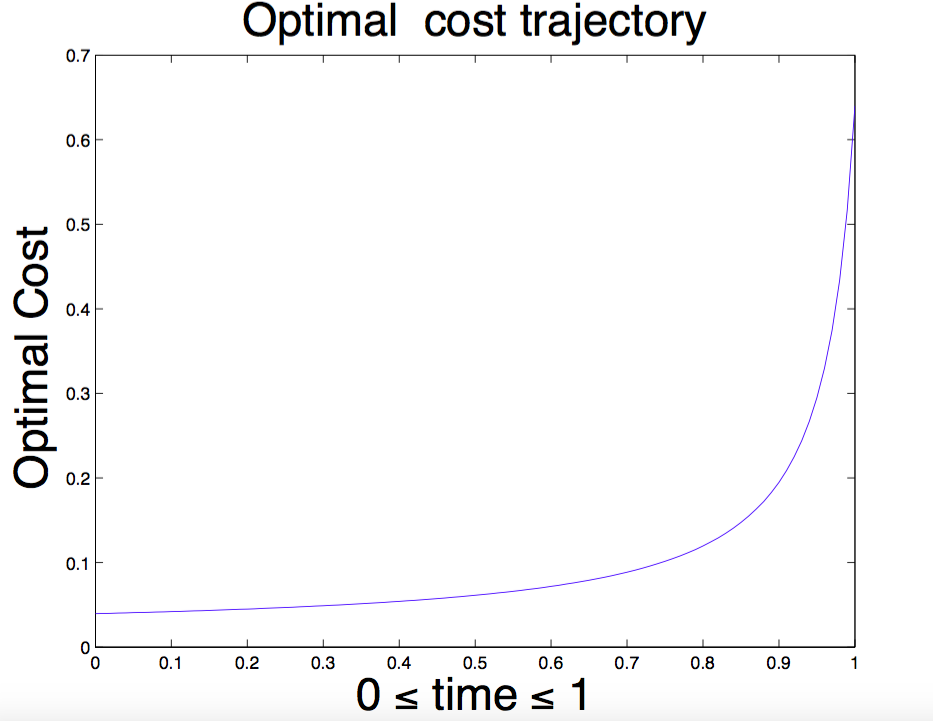}\\
  \caption{Optimal cost over time.  }\label{rnstateneu}
\end{minipage}
\end{figure}

\begin{figure}[htb]
\centering
\begin{minipage}{0.9\linewidth}
 \includegraphics[width=0.8\textwidth]{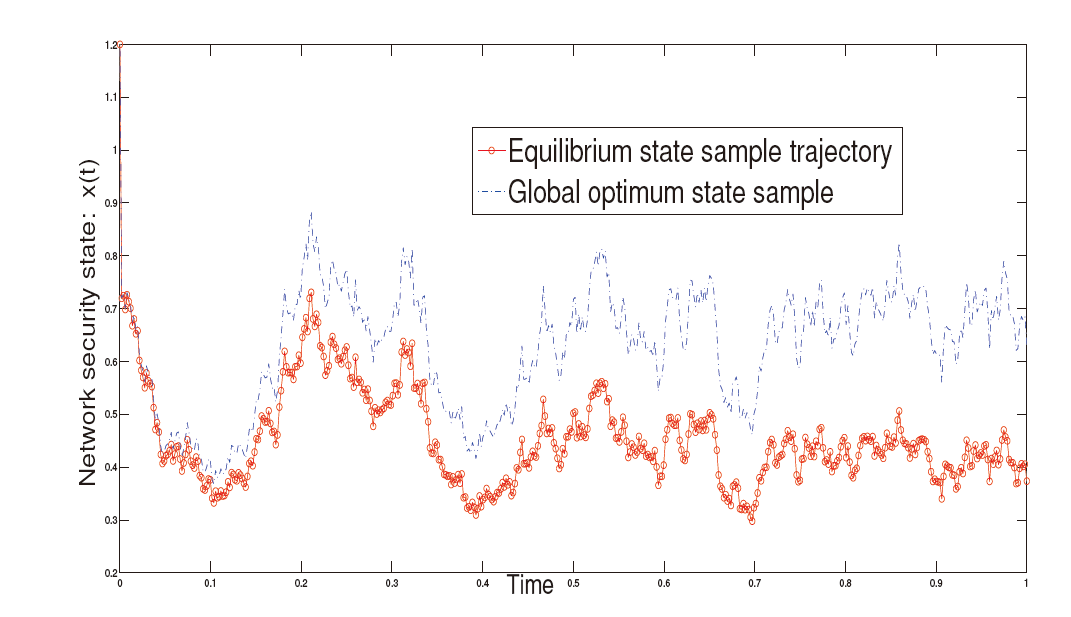}\\
  \caption{Optimal vs equilibrium state trajectory over time. The security level induced at equilibrium  state is lower than the one induced at full cooperation.  }\label{rnstateneu3}
\end{minipage}
\end{figure}
\end{example}

%%%%%%%%%%%%%%%%%%%%%%%%%%%%%%%%%%%%%%%%%%%%%%%%%%%%%%%%%%%%%%%%%

%%%%%%%%%%%%%%%%%%%%%%%%%%%%%%%%%%%%%%%%%%%%%%%%%%%%%%%%%%%%%%%%%%%

 Example \ref{example:security} can be used in the discrete-time mean-field-type game problem (\ref{smp0best}) associated with  (\ref{eq:discrete}). It corresponds to a variance reduction problem which is widely used in risk quantification. The following example solves a distributed  variance reduction problem in discrete time using MFTG.
\vspace{-2mm}
\begin{example}[Distributed Mean-Variance Paradigm, \cite{k102}] \label{variance2}
\emph{The best response problem  of agent $i$ is}
  \begin{eqnarray}
\label{Markowitz}   
\left\{ 
\begin{array}{lll} \vspace{-5mm}\\ \inf_{u_i\in \mathcal{U}_i}   \left\{ q_{iT}var(x_T)+(q_{iT}+\bar{q}_{iT})(E[x_T])^2\right. \\ \left. +\sum_{t=0}^{T-1}q_{it} var(x_t)+({q}_{it}+\bar{q}_{it}) (E[x_t])^2  + 
\sum_{t=0}^{T-1} r_{it}var(u_{it})+r_{it} (E u_{it})^2 \right\}
\\
\displaystyle{\mbox{ subject to }\ }\\
x_{t+1}=\left\{ax_t+\bar{a}Ex_t+\sum_{i=1}^n b_i u_{it}\right\}+\sigma W({t}),\ \\ x_0\sim \mathcal{L}(X_0),\  \ E[X_0]=m_0,  \\ 
 \end{array}
\right.
\end{eqnarray}
\emph{given the strategies $(u_j)_{j\neq i}$ of the other agents than $i$.}

\emph{Under the  assumption that for $t\in \{0,\ldots, T-1\},$ and
$ q_{jt}\geq 0,\ (q_{jt}+\bar{q}_{jt})\geq 0, \ r_{jt}>0,$  there exists a unique best-response  of agent $i$ and it is given by}
\begin{equation} \left\{\begin{array}{lll} 
{u}_{i,t}=  \eta_{it}(x_t-Ex_t)+\bar{\eta}_{it} Ex_t,  \\
 \eta_{it}=-\frac{[ ab_i {\beta}_{i,t+1}+b_i{\beta}_{i,t+1}\sum_{j\neq i }b_j\eta_{jt}]  }{{r}_{it}+b_i^2{\beta}_{i,t+1}},\\
\bar{\eta}_{it} =-\frac{b_i{\gamma}_{i,t+1}  (a+\bar{a}+\sum_{j\neq i}b_j\bar{\eta}_{j,t})}{ r_{it}+b_i^2{\gamma}_{i,t+1}},\\
{\beta}_{it}=q_{it} +{\beta}_{i,t+1}\{ a^2+2a\sum_{j\neq i }b_j\eta_{jt} +[ \sum_{j\neq i }b_j\eta_{jt} ]^2\}   -\frac{[ a b_i {\beta}_{i,t+1}+b_i{\beta}_{i,t+1}\sum_{j\neq i }b_j\eta_{jt}]^2  }{{r}_{it}+b_i^2{\beta}_{i,t+1}},  \\
{\beta}_{iT}= q_{iT}\geq 0\\
{\gamma}_{it}=({q}_{it}+\bar{q}_{it})+{\gamma}_{i,t+1}(a+\bar{a}+\sum_{j\neq i}b_j\bar{\eta}_{j,t})^2-\frac{(b_i{\gamma}_{i,t+1}  (a+\bar{a}+\sum_{j\neq i}b_j\bar{\eta}_{j,t}))^2}{ r_{it}+b_i^2{\gamma}_{i,t+1}},\\
{\gamma}_{iT}=  q_{iT}+\bar{q}_{iT}\geq 0\\
\end{array}\right. \end{equation}
\emph{and the best response cost of agent $i$ is}
$$
E[ L_i({u})]= E{\beta}_{i0}(x_0-Ex_0)^2+{\gamma}_{i0}(Ex_0)^2+ \sum_{t=0}^{T-1}{\beta}_{i,t+1}\sigma^2.$$
\end{example}

In both examples \ref{example:security} and \ref{variance2} the optimal strategy of agent $i$ is a feedback function of the state and the expected value of the state. This  structure is different than the one obtained in classical stochastic optimal control which are mean-field-free. The methodology used in standard stochastic  game problems do not apply directly to the mean-field-type game problems. These techniques need to be extended. This leads new optimality systems \cite{alainv1,pri6}. 
\section{ Engineering Applications}
\subsection{Civil Engineering}
This subsection discusses two applications of MFTG in civil engineering. 
\vspace{-2mm}
\begin{appli}[Road Traffic over Networks ]
\emph{The example below concerns transportation networks under dynamic flow and possible stochastic incidents on the lanes.
Consider a network ${\cal( V,L )}$, where
$\cal{V}$ is a finite set of nodes and ${\cal L \subseteq V
\times V } $ is a set of directed links. $n$ users share the network
$\cal (V,L)$. Let ${\cal R}$ be the set of
possible routes in the network. A user with a given source-destination pair arrives in the system at source node $s$ and
leaves it at the destination node $d$ after visiting a series of nodes and links, which
we refer to as a route or path. Denote by $c^{w}_i(x_t,u_{it}, m_t)$ the average $w-$weighted cost for
the path $u_{it}$ when $m_t$ fraction of users choose that path at time $t$ and $x_t$ is the incident state on the route.
The weight $w$ simply depicts that the effective cost is the weighted
sum of several costs depending on certain objectives. These metrics could be the delayed costs, queueing times, memory costs, etc and   can be
weighted by $w$ in the multi-objective case.
% Again, the weight $w$ could be
%different for different users due to their objectives. Henceforth,
%we omit $w$ and work with generic cost $c_i(x_t,u_{it},m_t)$ for
%simplicity of notation. We assume that the
%cost is non-decreasing in the  variable $m_t$ (congestion effect). 
We define two regimes for the traffic game: a finite regime
game with $n$ drivers denoted by $\mathcal{G}_{n}$ and an infinite regime game
denoted by $\mathcal{G}_{\infty}.$
The basic components of these games are $(\mathcal{N}, \mathcal{X},\mathcal{R}, I=\{x\}, c_i(x,.)).$
A pure strategy of  driver $i$ is a mapping from the
information set $I$ to a choice of a route that belongs to $\mathcal{R}.$ The
set of pure strategies of a  user  is $\mathcal{R}^{\mathcal{X}}.$ }

\emph{An action profile (route selection) $(u_1,\ldots,u_n)\in\mathcal{R}^{n}$ is an equilibrium of the finite mean-field-type  game if for every user $i$ the following holds:}
\[
c_i(x,u_{i},m(x,u_{i}))\leq c_i(x,u^{\prime
}_i,m(x,u^{\prime}_i)+\frac{1}{n}), \forall u^{\prime}_i\in\mathcal{R},
\]  \emph{for the realized state $x.$ }

\emph{The term $+\frac{1}{n} $ is  the contribution of the deviating user to the new route.
When $n$ is sufficiently large the state-dependent equilibrium notion becomes 
a population profile $m(x)=(m(x,u))_{u\in\mathcal{R}}$  such that for every user $i$}
\[
 m(x,u)>0 \Longrightarrow
c_i(x,u,m(x,{u}))\leq c_{i}(x,u^{\prime},m(x,u^{\prime})), %
\] \emph{for the realized state $x$ and for all  $ u^{\prime}\in\mathcal{R}.$ 
We refer to the equilibrium defined above as $0-$Nash equilibrium. Note that  the equilibrium profile depends on the realized state $x$. }

\begin{figure}[h]
\centering
\includegraphics[width=10cm, height=9cm]{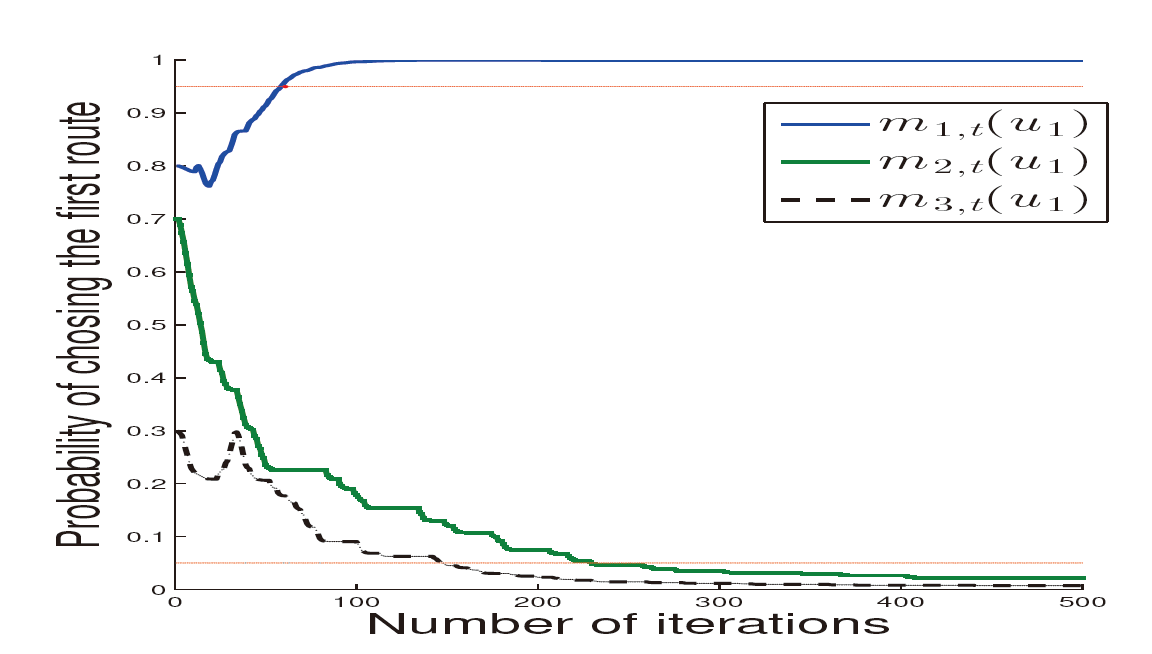}\\
\caption{Evolution of strategies of three agents over time. The imitative mean-field learning converges to a global optimum.} 
\label{fig:route1}
\end{figure}

\begin{figure}[h]
\centering
\includegraphics[width=10cm,height=7cm]{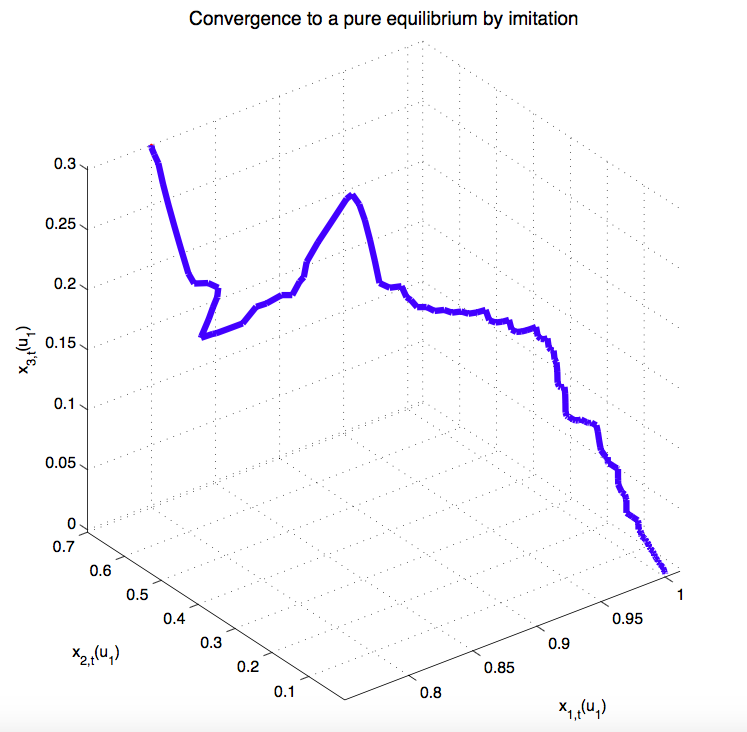}\\
\caption{The imitative mean-field learning converges to a global optimum.} \label{fig:route2}
\label{fig:route2}
\end{figure}

\emph{We now discuss the existence conditions. }

\emph{The equilibrium conditions can be rewritten in the form of variational inequalities:
for each state} $x,$ $(*)\ \sum_{u\in\mathcal{R}}[m(x,u)-y(x,u)] c(x,u,m(x,u))\leq 0,\ $ \emph{for all $y.$ }
\emph{ Hence, the existence of an equilibrium  is reduced to the existence of a solution to the variational inequality (*). By  the standard fixed-point arguments, we know from \cite{tembi09} that}
\emph{for each single state, such a population game has an equilibrium if the cost functions are continuous in the second variable $m$. Moreover, the equilibrium is unique under strict monotonicity conditions of the cost function $c_i(x,u,.).$ Note that uniqueness in $m$ does not mean uniqueness of the action profile $u$ since one can permute some of the commuters. We use  imitative learning in an information-theoretic view point. We introduce the cost of learning from strategy $m_{i,t-1}$ to $m_{i,t}$ as the relative entropy}
$d_{KL}(m_{i,t-1},m_{i,t}).$  %%%\footnote{$d_{KL}(x,y)=\sum_{y_i>0} y_i \log(\frac{y_i}{x_i}).$} 

\emph{Then, each user reacts by taking a myopic conjecture
given by $$\min_{m_{i,t}}\ \langle \hat{c}_{i,t}, m_{i,t}\rangle +\frac{1}{\beta_{i,t}}d_{KL}(m_{i,t-1},m_{i,t})$$ where $\hat{c}_{i,t}$ is the estimated cost vector,  $\beta_{i,t}$ is a positive parameter,
$d_{KL}$ is the relative entropy from }  $m_{i,t-1}$ to $m_{i,t}.$ 

%%%\footnote{$\beta_{i,t}$ is sometimes interpreted as temperature parameter or rationality level of driver $i.$ The
%%Boltzmann-Gibbs distribution can be obtained from noisy delay with $i.i.d$
%%measurement noise with Gumbel law. However, such a distribution is not observed in our measurement noise in practice}

\emph{$d_{KL}$ is not a distance (because it is not symmetric) but it is positive and can be seen as a cost to move from  $m_{i,t-1}$ to $m_{i,t}.$ We use the convexity property of the relative entropy to compute the strategy that minimizes the perturbed  expected cost.}
%\begin{proposition}
%The minimizer of $\ \langle \hat{c}_{i,t}, m_{i,t}\rangle +\frac{1}{\beta_{i,t}}d_{KL}(m_{i,t-1},m_{i,t})$ is the imitative Boltzmann-Gibbs strategy given by
%$$m_{i,t}(u):=\frac{ m_{i,t-1}(u)
%e^{-\beta_{i,t} \hat{c}_{i,t-1}(u)}}{\sum_{u^{\prime}\in\mathcal{R}}
%m_{i,t-1}(u^{\prime})e^{-\beta_{i,t} \hat{c}_{i,t-1}(u^{\prime})}}$$
%\end{proposition}
%By direct computation, one obtains that the minimizer strategy can be written as   multiplicative weighted imitative Boltzmann-Gibbs strategy.
\vspace{-2mm}
\begin{proposition}
\emph{Let $\beta_{i,t}=\log(1+\nu_{i,t})$ for $\nu_{i,t}>0.$ Then,
the imitative Boltzmann-Gibbs strategy is the minimizer of the above problem which becomes a multiplicative weighted imitative strategy:}
$$m_{i,t}(u):=\frac{ m_{i,t-1}(u)
(1+\nu_{i,t})^{ -\hat{c}_{i,t-1}(u)}}{\sum_{u^{\prime}\in\mathcal{R}}
m_{i,t-1}(u^{\prime})(1+\nu_{i,t})^{- \hat{c}_{i,t-1}(u^{\prime})}}.$$
\end{proposition}

\emph{The advantage of the imitative strategy is that it makes sense not only in small learning rate   but also in high learning rate. 
When the learning rate is large, the trajectory gets closer to  the best reply dynamics and for small learning it leads to the replicator dynamics \cite{taylorjonker1978}. One useful interpretation of the imitative strategy is the following: Consider a bounded rationality setup where the parameter $\nu_{i,t}$ is the rationality level of user $i.$ Then, a large value of $\nu_{i,t}$ means a very high rationality level for user $i,$  hence user $i$ will use an almost  ``best reply'' strategy. Small value of $\nu_{i,t}$ means that user $i$ is of a low rationality level  and is described by the replicator equation. It is interesting to see that both behaviors can be captured by the same imitative  mean-field learning.
 Note that the logit  (or Boltzmann-Gibbs) learning does not cover the low rationality level case.}

\vspace{-2mm}
\begin{proposition}
\emph{As $\nu_{i,t}$ goes to zero, the trajectory of the multiplicative weighted imitative strategy is approximated by the replicator equation of the estimated delays} $$
\dot{m}_{i,t}(u)=m_{i,t}(u)\left[ -\hat{c}_{i,t}(u)+\sum_{u^{\prime}} m_{i,t}(u^{\prime})\hat{c}_{i,t}(u^{\prime})  \right].
$$
\end{proposition}

\emph{For one commuter case, the solution of the replicator equation yields}
$$
m_{i,t}(u)= \frac{m_{i,0}(u)e^{-t .\frac{1}{t}\int_0^t \hat{c}_{i,t'}(u)\ dt'}}{\sum_{u^{\prime}}
m_{i,0}(u^{\prime})e^{-t .\frac{1}{t}\int_0^t \hat{c}_{i,t'}(u^{\prime})\ dt'}}
$$

\emph{The solution is}
$$
m_{i,t}(u)= \frac{m_{i,0}(u)e^{-t  \bar{c}_{i}(u)}}{\sum_{u^{\prime}}
m_{i,0}(u^{\prime})e^{-t  \bar{c}_{i}(u^{\prime})}}.
$$

\emph{Clearly the time-average trajectory  based on average payoff and smooth best reply dynamics are closely related with parameter $\beta_{i,t}=t.$
Each driver knows the current state  and employs the
learning pattern. Each driver
tries to exploit the information on the current state and build a strategy
based on the observation of the vector of realized delays over all the routes
at the previous steps. Then the Folk theorem for evolutionary game dynamics states:}
\begin{itemize}\item \emph{When starting from an interior mixed strategy, the replicator  equation converges to one of the equilibria.}
\item \emph{All the faces of the multi-simplex are forward invariant. In particular, the pure strategies are steady states of the imitative dynamics.}
    \item \emph{The set of global optima belongs to the set of steady states of the imitative dynamics.}
\end{itemize}

\emph{ The strategy-learning of   user $i$ is
given by}
{\small \begin{equation}
\label{tecost1}\mathcal{L}_{i}^1(x_{t}):\ \ \ m_{i,t}(x_{t},u):=\frac{
m_{i,t-1}(x_{t},u)(1+\nu_{i,t})^{- c_{i,t-1}(x_{t},u)}}{\sum_{u^{\prime}\in\mathcal{R}} m_{i,t-1}(x_{t},u^{\prime})
(1+\nu_{i,t})^{- c_{i,t-1}(x_{t},u^{\prime})}}%
\end{equation}

\begin{equation}
\label{tecost2}\mathcal{L}_{i}^{2}(x_{t}):\ \ \ m_{i,t}(x_{t},u):=\frac{ m_{i,t-1}(x_{t},u)
(1+\nu_{i,t})^{ -\bar{c}_{i,t-1}(x_{t},u)}}{\sum_{u^{\prime}} m_{i,t-1}(x_{t},u^{\prime}) (1+\nu_{i,t})^{-
\bar{c}_{i,t-1}(x_{t},u^{\prime})}},
\end{equation}
}
\emph{where $\bar{c}_{i,t}(x,u)$ is the time-average delay (up to $t$) in route $u$
and state $x.$}

\emph{The imitative mean-field learning above can be used to solve a long-term mean-field game problem. We observe in Figures \ref{fig:route1}- \ref{fig:route2} that the imitative learning converges to one of the global optima. However, the exploration space grows in complexity. We explain how to overcome to this issue using mean-field learning based on particle swarm optimization (PSO). In it each user has a population of particles (multi-swarm). The particles within the same population (coalition) may pool their effort to learn faster and exploit better the available information.}

\emph{The next example concerns multi-level building evacuation \cite{aims,burj2,burj3,burj4} using constrained mean-field games.}
\end{appli}

\vspace{-5mm}
\begin{appli}[Multi-level building evacuation]
\emph{A typical  mean-field game model assumes that agents have unconstrained state dynamics. This has been, for example, the case with most of the existing mean-field models developed in the last  
three decades.  
Such models may not however be useful in practice, for example in a context of building evacuation.   Evacuation strategies and values are designed 
 using constrained mean-field-type game theory.}  
%%%%%%%%%%%%%%%%%%%%%%%%%%%%%%%%%%%%%%%%%%%%%

\begin{figure*}[ht]
\centering
\includegraphics[width=5cm,height=4cm]{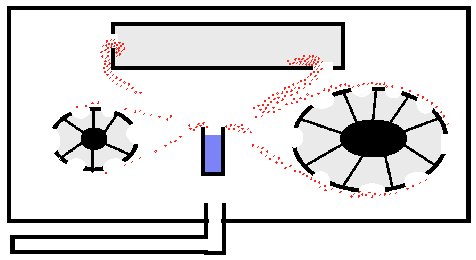}\
\includegraphics[width=5cm,height=4cm]{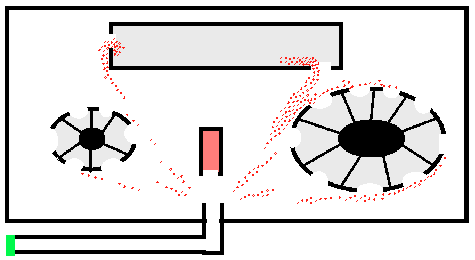}
\caption{Spatial distribution of agents at time $t= 5.$ Agents are represented by small circles in the map. Agent in the higher floors will be
evacuated using the stair (blue rectangle) on floor 2. There is one exit door in the ground floor. The exit door is in green-color code in the ground floor. Each agent chooses the shortest  and less congested path and decreases its velocity according to its own congestion measure.\label{fig:minipage2}}
\end{figure*}

\begin{figure*}[ht]
\centering%
\subfigure 
{ \includegraphics[width=8cm,height=5cm]{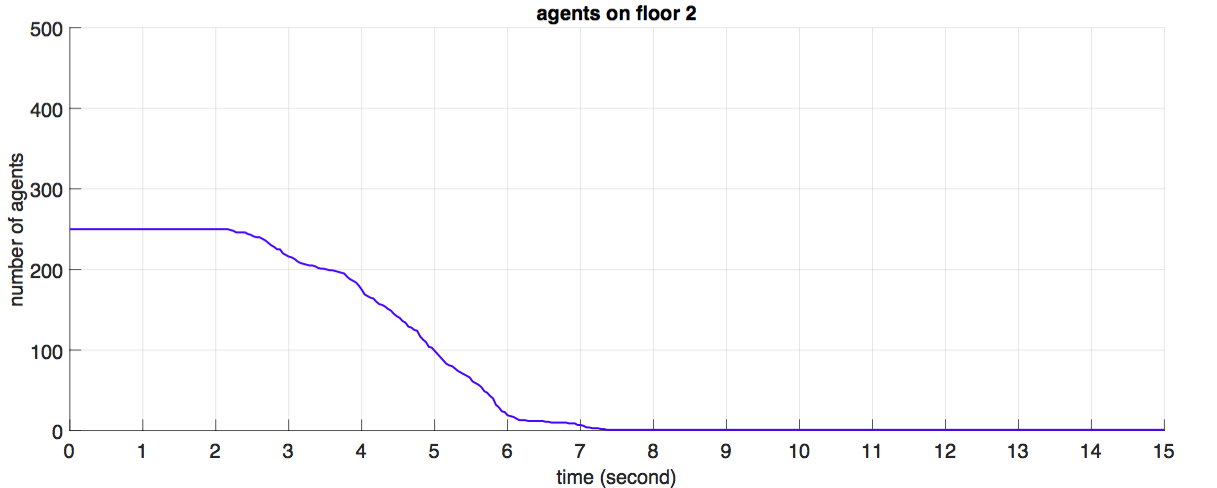} } 
\subfigure 
{ \includegraphics[width=8cm,height=5cm]{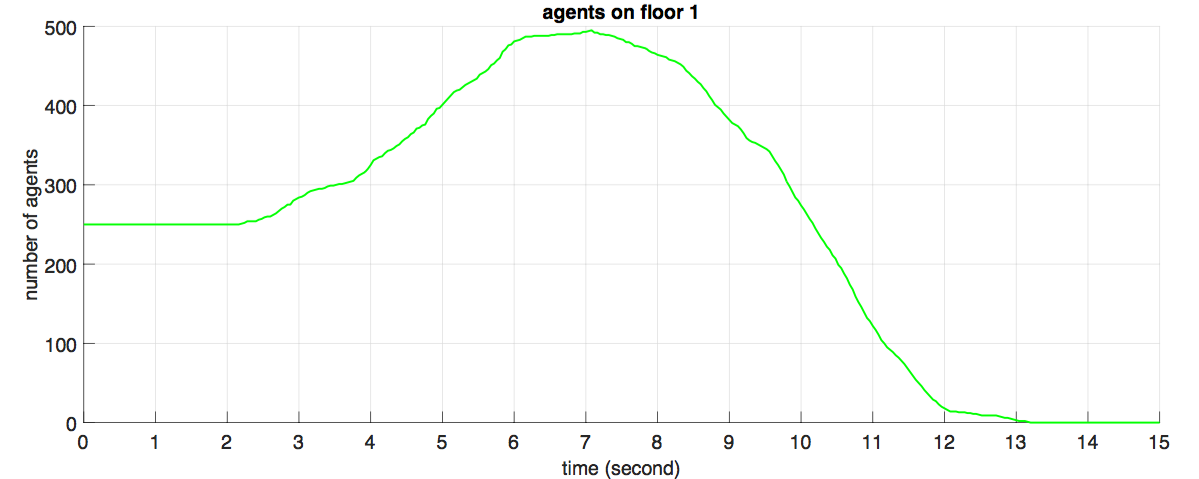}  }\\
\subfigure
{  \includegraphics[width=8cm,height=5cm]{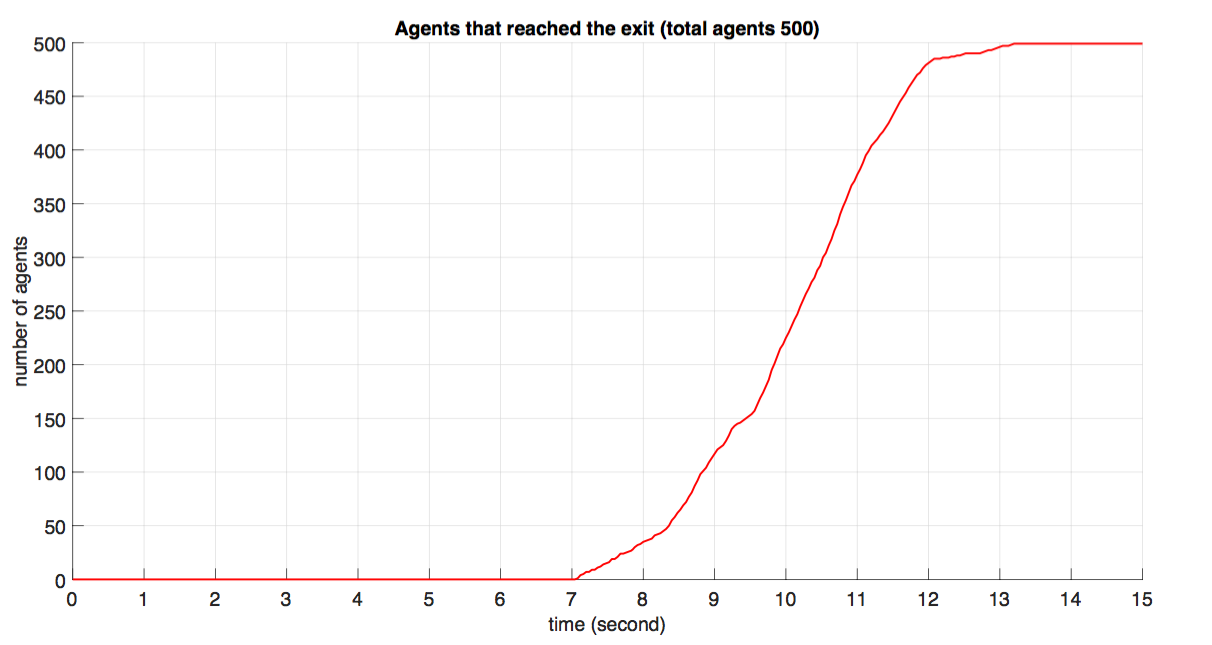}  }
\caption{{The two upper figures plot the evolution of the number of remaining agents in the building. The number of agents in ground floor starts increasing because the flow is coming from first floor until certain time threshold and then decrease when agents start to exit.  The lower Figure plots
the evolution of the number of agents who have been evacuated safely. The plot has a typical shape of a cumulative distribution function.} \label{fig:t123}}
\end{figure*}
%%%%%%%%%%%%%%%%%%%%%%%%%%%%%%%%%%%%%%%%%%%%%
\emph{Particle-based pedestrian models have been  studied in  \cite{ref10,ref12}. Continuum approximation of theoretical models have been proposed in \cite{ref9,ref10,ref11,ref12,ref13,cedric}. Recent mean-field studies on  crowd and pedestrian flows include \cite{dogbe,wolfram,marie1,marie2,diogoobstacle}.   {\color{black} Below a mean-field game for  multi-level building evacuation is presented. Consider a building  with multiple floors and resolutions represented by a compact domain $D$ in the $m-$dimensional Euclidean space $\mathbb{R}^m.$ The number of floors is $K.$ The domain at floor $k$ is denoted as $D_k.$ For $1<k<K,$ the floor $k$ is connected to the higher floor  $k+1$ using the intermediary domain $I^{+}_{k}$ but also the lower floor $k-1$ using $I^{-}_{k}.$ The sets $I_k$ can be elevator zones or stairs.
$n\geq 2$ agents are distributed in a multi-level multi-resolution building with stairs, exit doors, sky-bridges. Each agent knows her current location in the building.
The state/location $x_i$ of an agent  $i$ changes depending on her control action $u_i$. The agent is interested in a safe evacuation from the building. This means that she is interested in the minimal exit time that avoid huge crowd around her.  The problem of the agent $i$ is equivalent to }
$$
\left\{
\begin{array}{ll}
  &  \inf_{u_i} \    c_3(x_i(T),G_n(x_i(T))) +\int_0^T c_1(G_n(x_i(t))) \|u_i(t)\|^2+ c_2(G_n(x_i(t))) \ dt,    \\
& \dot{x}_i=u_i\in  \mathbb{R}^3 ,\    0<t<T \\  
&  x_i(t) \in D \subset \mathbb{R}^3,  \\  
 & \mbox{Boundary constraints:} \\  &  u_{i | \ \partial D }=0,\    u_{i | \ \mathrm{ExitDoor} }=k\neq 0  \in \mathbb{R}^3
\end{array}
\right.
$$ where $c_i$  is a positive increasing function, with $c_2(0)=0.$  $T>0$ is the  exit time at one of the exits. The final exit cost is represented by $c_3$ which can be written as $\tilde{c}_3+\tilde{h}(x)$ where $\tilde{c}_3>0$ captures the initial response time of an agent (without congestion around), 
$$G_n(x_i(t))= \frac{1}{vol(B(x_i(t),\epsilon))}\ \frac{\sum_{j\neq i} \ind_{\{d(x_j(t), x_i(t))\leq \epsilon\}}}{n-1},$$ represents the number of the agents around the position $x_i$ except $i$ within a distance less than $\epsilon>0,$ $vol(B)$ is the $m$-dimensional volume of the ball $B(x_i(t),\epsilon)$ which does not  depend on $x_i(t),$ due to translation invariance of the volume measure.  When the number of agents grows, one obtains a mean-field game with several interacting agents. The state dynamics must satisfy the constraint $x_i(t)\in D$ at any time $t$ before the exit.}
\emph{ The  non-optimized Hamiltonian in macroscopic setting as $$H^0(x,u,G,p)=-c_1( G(x) )\|u\|^2 -c_2(G(x))+p.u,$$ where $p$ is the adjoint variable. 
 The Pontryagin maximum principle yields}
$$\begin{array}{c}
 \dot{p}= -H_{x}^0, \\ 
 p(T)=-g_{x}(x(T)), 
 \\
 \dot{x}=  \frac{p}{2 c_1(G(x)) },\   0< t\leq T \\ x(0)\in D. 
\end{array}$$

\emph{The Hamiltonian $H^0 ( ., ., G, p(t))$ is concave in $(x, u) $ for almost everywhere (a.e.)  $t \in [0, T ].$
Then, for   convex function $c_3,$  $u^*$ is an optimal  response if $H^0(x^*(t), u^*,G^*,p^*(t))= \max_{u}H^0(x^*(t), u,G^*,p^*(t)).$
The (optimized) Hamiltonian as}
$$
H(x,p,G)=\sup_{u} \{ - c_1(G(x))\| u\|^2- c_2(G(x)) +p u\}.
$$
\emph{ The Hamiltonian can 
 be computed as
$
H(x,p,G)=   \frac{\|p\|^2}{ 4 c_1( G(x))}-c_2(G(x)),
$
and the optimal strategy is in (own)state-and-mean-field feedback form: $u^* =\frac{p}{2 c_1(G(x)) }=H_p(x,p,G(x)), $ to be projected to the tangent space.
The dynamic programming principle  leads to the following optimality system:}
$$ \left\{ \begin{array}{c}
v_t+H(x,v_x,G(x))=0, \ \mbox{ on}\   (0,T)\times {D}\\
v(T,x)=-g(x),   \mbox{ on}\   D\\
\rho_t+div_x(\rho H_p)=0,\ \ 
\rho_0(.)  \mbox{ on}\   D  \subset \mathbb{R}^3 \\
 u=0, \ \ y=0 \mbox{ on}\  \partial D \\ u=k, \  \mbox{ at}\  \mbox{exits} \ \end{array}
\right.
$$

\emph{The development of numerical result,  simulation and a validation framework  can be found in  \cite{aims,burj2,burj3,burj4}. Figures   \ref{fig:minipage2} and  \ref{fig:t123} show the application to a two floors  building where 500 agents are spatially distributed.}
 
\end{appli}
\vspace{-3mm}
Next, two  applications of MFTGs in electrical engineering are presented.
\subsection{Electrical Engineering}
\begin{figure}[htbp]
\begin{center}
    \includegraphics[width=12cm, height=8cm]{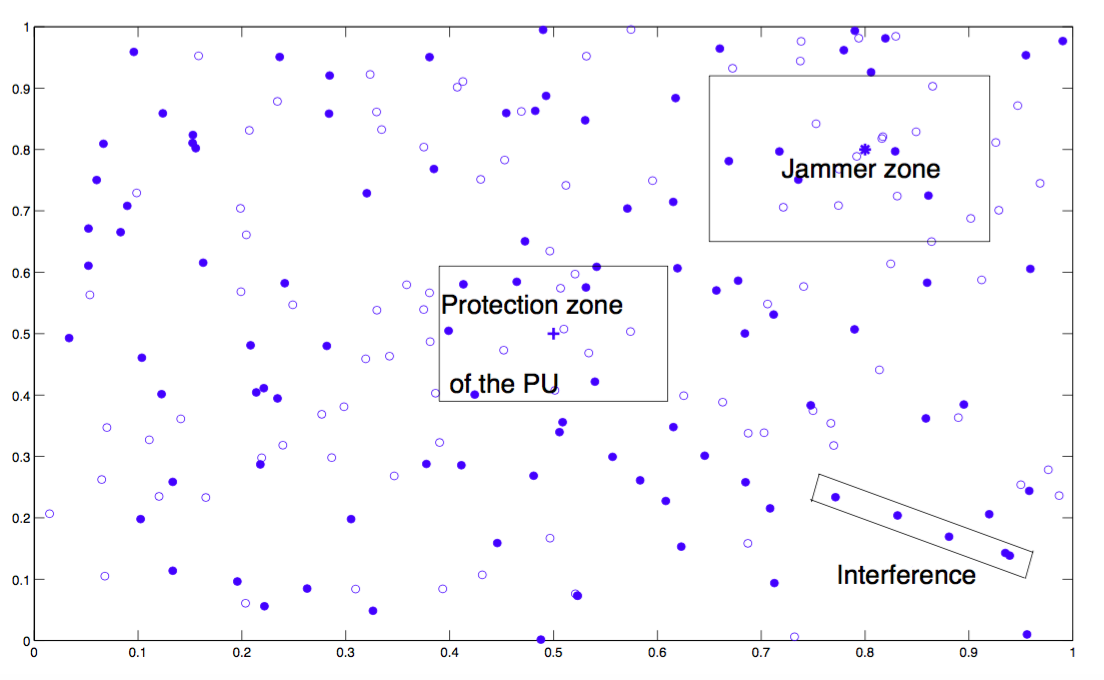}\\
  \caption{A typical large-scale network with regular nodes, relay nodes, primary users and  jammers. The star sign represents a Jammer. The blue
  nodes are active secondary nodes, the nodes in circle are inactive secondary nodes, and the plus sign represents a primary user zone of transmission using MIMO millimeter wave wireless communication.}
  \label{wave}
  \end{center}
\end{figure}

\begin{appli}[Millimeter Wave Wireless Communication]
\emph{Millimeter wave (mmWave) frequencies, roughly between 30 and 300 GHz, offer a new frontier for wireless
networks. The vast available bandwidths in these frequencies combined with large numbers of spatial degrees
of freedom offer the potential for orders of magnitude increases in capacity relative to current networks
and have thus attracted considerable attention for next generation 5G communication systems. 
However,
sharing of the spectrum and the available infrastructure will be essential for fully achieving the potential
of these bands. Unfortunately, rapidly changing network dynamics make it difficult to optimize resource
sharing mechanisms for mmWave networks. MIMO mmWave wireless networks will rely extensively on highly directional transmissions, where both users, relays
and base stations transmit in narrow, high-gain beams through electronically steerable antennas. While
directional transmissions can improve signal range and provide greater degrees of freedom through spatial
multiplexing, they also significantly complicate spectrum sharing. Nodes that share the spectrum must not
only detect one another, but also search over a potentially large angular space to properly steer the beams
and reduce interference. 
Power allocation, angle optimization and channel selection algorithms should consider the possible interference field and reduce it by adjusting  the angles. This can facilitate rapid
directional discovery in a dynamic and mobile environment as in Figure \ref{wave}. 
Sometimes jammers and malicious are involved in the interactions. 
Beams adjustment and Interference coordination are  central problem
 for users within the same network, or between users in different networks sharing the same
spectrum. When multiple operators own separate core network and radio access network (RAN) nodes such
as base stations and relays, but only loosely coordinate via wireless signaling, it is essential to use incentive
mechanisms for better coordination to exploit the available resources. Cost sharing and pricing mechanisms
capture some of the fundamental properties that arise when sharing resources among multiple operators. It
can also be used in the uplink case, where users can select their preferred services and network provides
and have to find tradeoffs between quality-of-experience (QoE) and cost (price). }

\emph{As an illustrative example, a  particle swarm learning mechanism, which is mean-field dynamics, in which the particles adapt the parameters such as angle and power  is used to  improve users' quality-of-experience. Here the key mean-field terms are the distributions of remaining energy, distribution of transmitter-receiver pairs and   the  sectorized  interference field (per angle). Since users are carrying  smartphones with limited power consumption, it is crucial to examine the remaining energy level. As in \cite{wireless4} 
the energy dynamic can be written as 
$$de= -u dt +v dt+\sigma dW,$$ subject to $e(t)\geq 0,$  $e(0)=e_0,$ and $u(.)\geq 0$ is the transmission power and $v(.)$ is the energy harvesting rate (for example, with distributed renewable energy sources).} 
\vspace{-3mm}
\begin{proposition} 
\emph{The marginal distribution $m^e(t,e)$ of remaining energy solves the Fokker-Planck-Kolmogorov equation:
$$ \partial_tm^e+\partial_e[(-u+v)m^e]-\frac{\sigma^2}{2}\partial_{ee}m^e =0,
$$
in a distribution sense. The first moment dynamics yields
$
\frac{d}{dt}\bar{e}=-\bar{u}+\bar{v},
$
where $\bar{e}(t)=\mathbb{E}[e(t)], \bar{u}(t)=\mathbb{E}[u(t)], \bar{v}(t)=\mathbb{E}[v(t)]$ denotes the expected value of $e(t)$, $u(t), v(t)$ respectively.}
\end{proposition}
\vspace{-3mm}
\emph{Users move according to a mobility dynamics (which may not be stationary).   The channel state can be modeled, for example using a matrix valued Ornstein-Uhlenbeck process 
$dH_j= \Gamma_j[\hat{H}_j-H_j] dt+dW_j$ where $\Gamma_j, \hat{H}_j$ are matrices with compatible dimensions of antennas at source and destination. }
\vspace{-3mm}
 \begin{proposition} 
\emph{The marginal distribution $m^{H_j}(t,H_j)$ of channel state of user $j$  solves the Fokker-Planck-Kolmogorov equation:
$$ \partial_t m^{H_j}+\mbox{div}_{H_j}[(\Gamma_j(\hat{H}_j-H_j))m^{H_j}]-\frac{1}{2}trace[\partial_{H_jH_j}m^{H_j}] =0,
$$
in a weak sense. The first moment dynamics of  user $j$ yields
$$
\frac{d}{dt}\bar{H}_j=\Gamma_j (\hat{H}_j-\bar{H}_j),
$$
which decays exponentially to $\hat{H}_j$ as $t$ increases.}
\end{proposition}

\vspace{-3mm}
\emph{The (unnormalized) distribution of the triplet (position, energy, channel)    of the population at time (or period) $t$ is $\nu(t,e,x,H)=\sum_{j=1}^n \delta_{\{e_j(t),x_j(t), H_j(t)\}},$ and the one within a beam $A(s,d)$ with direction $s-d$ is 
$$\tilde{\nu}(t,e,x,H,s,d)=\sum_{j=1}^n \delta_{\{e_j(t),x_j(t), H_j(t)\}}\ind_{x_j(t)\in A(s,d)}.$$
The sectorized interference field is $I(t,x(t),d)=\int_{(\bar{x}, \bar{H},\bar{u})}\phi(\bar{x}-x(t), \bar{H}, \bar{u}))\tilde{\nu}(t,E,\bar{x},\bar{H},x(t),d).$ 
Compared to other wireless technologies, mmWave may generate less interference because of reduced and optimized angles. However, interference  may still occur when several users and blocking objects fall within the same angle as depicted in Figure \ref{wave}.  The success probability  $\mathbb{P}(\mbox{SINR}_i \geq \beta_i)$ from position $x_i(t)$ to destination $d_i$ for both LoS and non-LoS can then be derived. The  quality-of-experience of users  can be termed as function of the sectorized interference field, satisfaction level and user-centric subjective measures such as MOS (mean opinion score) values.  }
\end{appli}

\vspace{-3mm}
\begin{appli}[Distributed Power Networks (DIPONET)]
\emph{ Distributed power is a power generated at or near the point of use. This includes technologies that supply both electric power and mechanical power. The rise of distributed power is also being driven by the ability of distributed power systems to overcome the energy need constraints, and transmission and distribution lines. 
Mean-field game theoretic applications to power grid can be found in \cite{pe1,pe2,pe3,pe4,pe5,pe6,pe7,pe8,fe8, fe9, fe10,ccdc16one,ye}. A prosumer (producer-consumer) is a user that not only consumes electricity, but can also produce and store electricity.  Based on forecasted demand, each operator determines its production quantity, its mismatch cost, and engages an auction mechanism to the prosumer market. The performance index is $ L_j(s_j,e_j)= l_{jT}(e(T)) +\int_0^T l_{j}(D_j(t)-S_j(t)) + \frac{\rho}{2}\sum_{k}s^2_{jk}(t)\ dt.$ Each producer aims to find the optimal production strategies:}
$$\left\{
\begin{array}{l}
 \ \inf_{s_j,e_j} L_j(s_j, e_j, T)\\
 \frac{d}{dt}e_{jk}(t)=c_{jk}(t) - s_{jk}(t)\\
 c_{jk}(t)\geq 0,\  s_{jk}(t)\in [0, \bar{s}_{jk}],\ \forall j, k, t\\
 s_{jk}(w)=0\ \mbox{if \ w is a starting time of a maintenance period. }\\
 e_{j,k}(0)\  \mbox{given}.
\end{array}
\right.
$$ \emph{where  $D_j(t)$ is a demand at time $t$,    $l_j(D_j(t)-S(t))$ denotes the instant loss  where $S(t)=S_{producer}(t)+S_{prosumer}(t),$  $S_{producer}(t)=\sum_{j=1}^n s_j(t)=\sum_{j=1}^n \sum_{k=1}^{K_ j}s_{j,k}(t)\ ,$  where $s_{j,k}(t)$ is the production rate of plant/generator $k$  of $j$ at time $t.$ $K_j$ total number of power plants of $j.$ The loss $l_j$ is assumed to be strictly convex. The stock of energy at time $t$ is given by the classical motion $\frac{d}{dt}e_{jk}$ where $c_{jk}(t)$ is the maintenance cost of plant/generator $k$ of $j$ when it is in the maintenance phase. 
The optimality equation of the problem is given by  Hamilton-Jacobi-Bellman:}
  \begin{equation}  \label{eqH}
    \begin{cases}
      \partial_{t} v_j(t,e_j) + H_j(D_j(t),  \partial_{e_j}v_j(t,e_j)  ) = 0 , \,\,   t < T   \\
      v_j(T,e_j) = l_{jT}(e_j), \,\,   
    \end{cases}
  \end{equation}  \emph{where $H_j$ is the Hamiltonian function is}

 \begin{eqnarray}H_j(D_j,y_j) = \inf_{s_j}[l_{j}( D_j-S_j) + \frac{\rho}{2}\sum_{k}s^2_{jk}+\sum_{k}(c_{jk}-s_{jk})y_{jk}]\
  \end{eqnarray}
\emph{ The first order interior optimality condition yields $- l'_{j}( D_j-S_j)-y_{jk}+\rho s_{jk}=0.$ 
 By summing over $k$ one gets an equation for the total production quantity $S_j^*$ solves
 $- K_j l'_{j}( D_j-S_j) -\sum_{k=1}^{K_j}y_{jk}+\rho S_j=0.$
 The optimal supply of power plant $k$ is 
 $
 s_{jk}^*= \min(\bar{s}_{jk}, \frac{l'_{j}( D_j-S^*_j)+y_{jk}}{\rho}).
 $
   The solution  of  partial differential equation (\ref{eqH}) can be explicitly obtained and it is  given by the Hopf-Lax formula: }
 \begin{equation} \label{thh}v_{j}(t,e_j)=\inf_{y\in  \mathbb{R}^{K_j}} \big\{ l_{jT}(y)+ (T-t)H^*_j\big (D_j,\frac {e_j-y}{T-t}\big)\big \},\end{equation}
\emph{ where $H^*_j$ is the Legendre transformation of  $H_j,$  and is given by}
 $$
 H^*_j(D_j, a)=l_{j}\left( D_j-\frac{1}{\rho}\sum_{k}a_{jk}-\frac{l'_{j}( D_j-S^*_j)}{\rho}\right) $$ $$+ \frac{\rho}{2}\sum_{k}a^2_{jk}+\sum_{k}c_{jk}a_{jk}.
 $$

\emph{ Note that  (\ref{thh}) provides an explicit solution to the Demand-Supply matching problem between power plants of  prosumer $j$  and this holds for arbitrary number of prosumers and power stations. }

   \begin{figure}[htbp]
\begin{center}
\includegraphics[width=15cm,height=9cm]{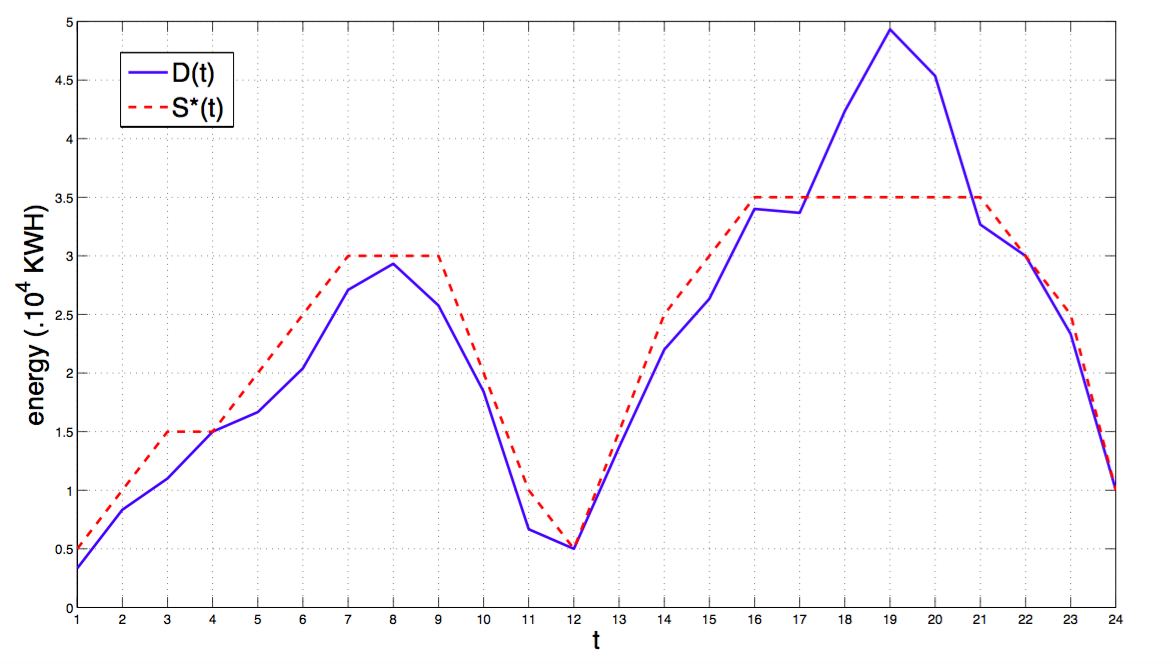}
\caption{{Optimal supply $S^*_j$ of producer $j$  obtained by means of inf-convolution of the Bellman operator}}
\label{defaultv03}
\end{center}
\end{figure}

 \begin{figure}[htbp]
\begin{center}
\includegraphics[width=15cm,height=8cm]{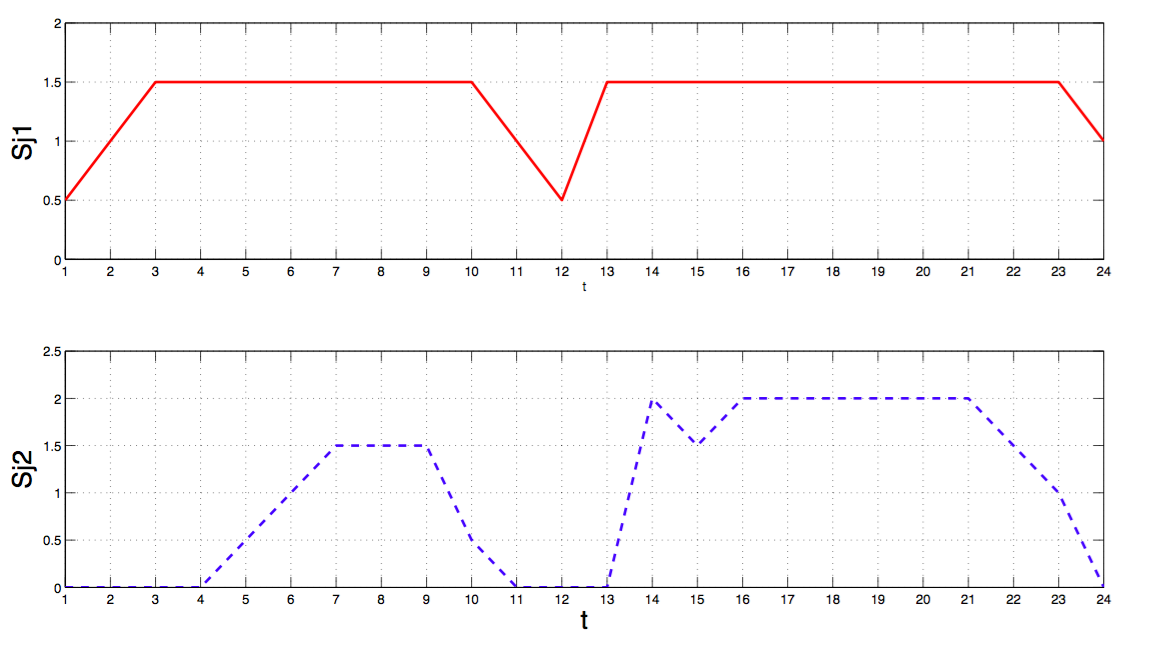}
\caption{{Optimal Allocation  $\sum_{k}s_{jk}(t)=S^*_j(t)$ between the two power stations of producer $j$ at time period $t$}}
\label{defaultv04}
\end{center}
\end{figure}

\emph{The mean-field equilibrium is obtained as fixed-point equation involving $S^*$ and $D^*.$ When $ l'_{j}$ is continuous and preserves the production domain $[0,\bar{s}]$ one can guarantee the existence of such a solution by using Brouwer fixed-point theorem. One can use higher order fast mean-field learning to learn and compute of such a mean-field equilibrium.
Figure \ref{defaultv03} illustrates the optimal supply based on an estimated demand curve. Figure \ref{defaultv04} represents an allocation of the producer with two power stations.} 
\end{appli}

\subsection{Computer Engineering}
This section provides applications of MFTG in computer engineering. It starts with an application of MFTG with number finite state-actions and  then focuses on continuous state-action spaces.

\vspace{-3mm}
\begin{appli}[Virus Spread over Networks]
\emph{We study a malware propagation over computer  networks where  the nodes interact through network-based opportunistic meetings (see Figure \ref{fig:rep0} and Table \ref{trans}). The security level of network is measured as a function of some key control parameters: acceptance/rejection of a meeting, opening/not opening a suspicious e-mail,  file or packet. We model the propagation of the virus in network as a sort of epidemic process on a random graph of opportunistic connections \cite{malware}. 
A computer/node  can  randomly get  online an infected or non infected data from other computers.}

\begin{figure}[htb] 
\begin{center}
\includegraphics[width=0.8\textwidth]{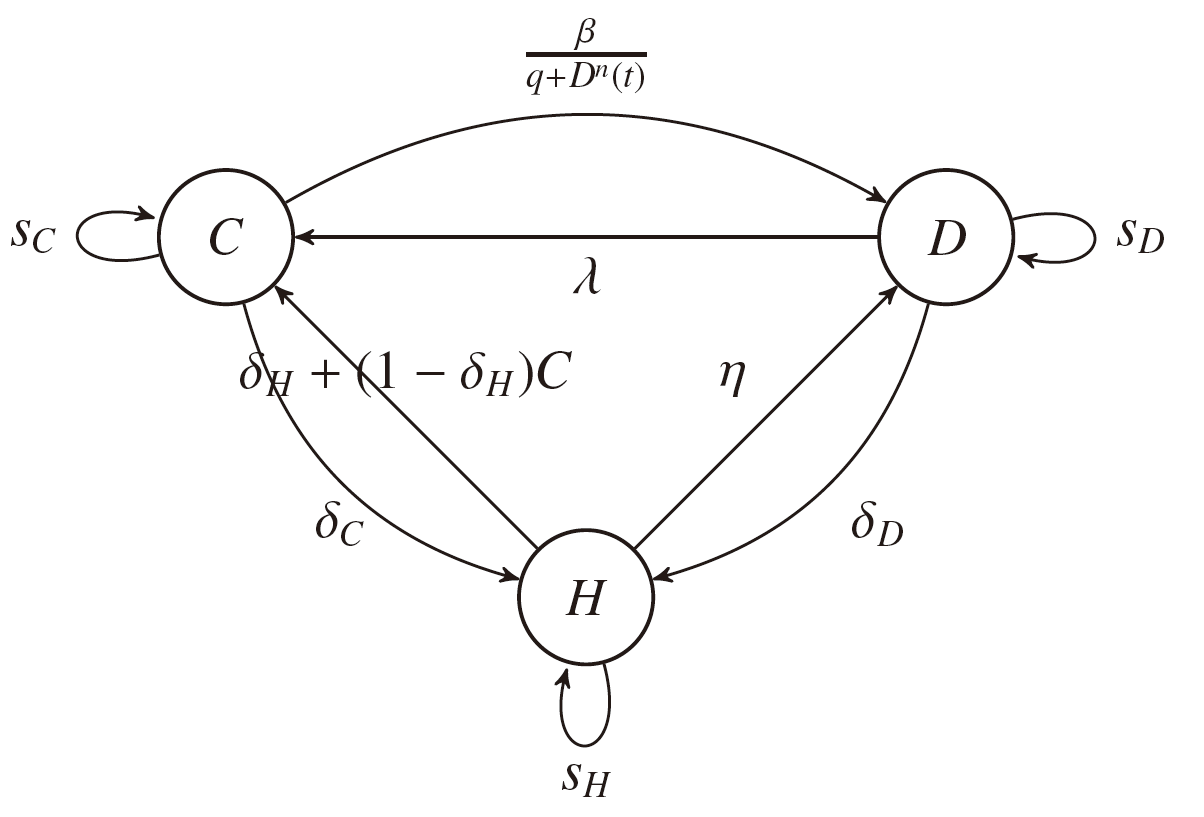}
\end{center}
\caption{ {Markov chain representation: the parameters $s_i$ are the complement of the other transitions.}}
\label{fig:rep0}
\end{figure}

\emph{An infected computer can be in two states: dormant or  fully infected. The non-infected computers are susceptible to be  approached  by virus coming from infected ones. The possible states are therefore denoted as Dormant (D), Infected/Corrupt(C) and Susceptible/Honest (H). The set of types is 1 or 2, also denoted generically as $\theta,\theta'.$ For each type the state may be different except for honest state where it is considered as honest in both regimes of the network.
The network size is  $n\geq 1.$ The repartition of the nodes at time step $t$ is denoted as  $n=D_{\theta}(t)+D_{\theta'}(t)+C_{\theta}(t)+C_{\theta'}(t)+H(t).$ }

\color{black} \emph{The frequency of the states  $\theta$ is called occupancy measure of the population and is denoted as 
$M^n(t) = (D_{\theta}(t)/n,D_{\theta'}(t)/n,C_{\theta}(t)/n,C_{\theta'}(t)/n,H(t)/n) =: (D^n_{\theta}(t),D^n_{\theta'}(t),C^n_{\theta}(t),C^n_{\theta'}(t),H^n(t)).$ $M^n(.)$ is a random process and its limit measure corresponds to the mean field term. The goal is understand the impact of the control action on combatting virus spread, which is the minimization of proportion $O^n(t): = 1-H^n(t))$. 
The interaction is simulated using the following rules:}

\emph{Changes from Dormant states: A node in dormant state (transient) with type $\theta$ may become honest  with probability $\delta_{D}\in (0,1).$
 A dormant  with type $\theta$ may opportunistically meet another dormant of type $\theta'$, and both become active. This occurs with probability proportional to the frequency of other dormant agent at time $t$.
For type $\theta,$ the probability is $\lambda (D^n_{\theta'}(t)- \frac{1}{n}\ind_{\{ \theta=\theta'\}})$. Note that the dormant can decide to contact the other dormant or not, so there are two possible actions: $\{m, \bar{m}\}$ (to meet or not  to meet). Those events will be modeled with a Bernoulli random variable with success (meeting) probability $\delta_m$, which represents $u(m| D,{\theta}).$ }

\emph{Changes from Corrupt States: A corrupt node may become honest  with probability $\delta_{C}$.
 A corrupt node of type $\theta$ may become dormant with probability $\beta \frac{D^n_{\theta}(t)}{q_{\theta} + D^n_{\theta}(t)}$ at time $t$. Here is assumed that, at high concentrations of dormants, each corrupt node infects at most a certain maximum number of dormant nodes per time step. This reflects the fact a corrupt has a limitation in terms its power, domination and capabilities. The parameter $0 \leq \beta \leq 1$  can be interpreted as a  maximum contamination rate. The parameter $0 \leq q_{\theta} \leq 1$ is the dormant  node density at which the infection spread proceeds. }

\emph{	Changes from Susceptible/Honest states:	 An honest node may become infected with probability $\delta_{H}+(1-\delta_{H})C^n (t).$ 
	 An honest node may become dormant via two ways. First, $\delta_{Sm}$ is the probability of getting corrupt by the network representative node. In this case, the honest node can decide share or not, so there are two possible actions: $\{o, \bar{o}\}$. This case will be modeled using  a coin toss with probability $\delta_e\in (0,1).$ Second, $\eta (D^n_{\theta}(t)+D^n_{\theta'}(t))$ models the probability of meeting a dormant node. 
	Here $\eta \in (0,1).$ In this case, the dormant node can decide to contact the honest node or not, and it is modeled analogously to the other two cases.}

{\color{black} \emph{The  payoff function is the opposite of the infection level. Each transition described above has a certain  contribution to be infection level of the society, which could be 0 if no corrupt or dormant node become honest, $-1/n$ if there is a node which become honest and $+1/n $ if one node is corrupt  (D or C). In Table \ref{trans} are the transition probabilities, the contribution to $M^n(t+1)-M^n(t)$, the set of actions, and the contribution to information spread in the network.}}

\begin{table}[ht]%
\begin{center}
\caption{{Probabilities, effects $(D,C,H)$, actions and loss function.}} \label{trans}
\begin{tabular}{lcccc}
\hline
Case & Transition proba. $(\theta,\theta' \in \{1,2\})$. & $M^n_{\theta}(t+1)-M^n_{\theta}(t)$ & Actions & Propagation \\ \hline
 ${D}  \xrightarrow{\delta_{D}} {H}$ & $D^n_\theta(t) \delta_{D}$ & $(-1,0,1)/n$ & singleton set & $-1/n$ \\
 $2 D \xrightarrow{\lambda} 2 C$ & $D^n_{\theta}(t) \delta_m ^2 \lambda (D^n_{\theta}(t)- \frac{1}{n})$ & $(-2,2,0)/n$ & $\{m, \bar{m}\}$ & $0$ \\
%2' & $P^n_{\theta}(t) \delta_m ^2 \lambda (P^n_{\theta'}(t)- \frac{1}{n}\ind_{\{\theta=\theta'\}})\ind_{\{\theta=\theta'\}} $ & $(-2,2,0)/n$ & $\{m, \bar{m}\}$ & $0$ \\

 ${C}  \xrightarrow{\delta_{C}} {H}$ & $C^n_\theta(t) \delta_C$ & $(0,-1,1)/n$ & singleton set & $-1/n$ \\
 ${C}  \xrightarrow{\frac{\beta}{q_{\theta} + D^n_{\theta}(t)}} D$ & $C^n_\theta(t) \beta \frac{D^n_{\theta}(t)}{q_{\theta} + D^n_{\theta}(t)}$ & $(-1,1,0)/n$ & singleton set & $0$\\
 ${H}  \xrightarrow{\delta_{H}+(1-\delta_{H})C^n} C$  & $H^n(t) [\delta_{H}+(1-\delta_{H})C^n(t)]$ & $(0,1,-1)/n$ & singleton set & $1/n$ \\
 ${H}   \xrightarrow{\eta}D$ & $H^n(t)(\delta_e \delta_{Sm} + \delta_m \eta  D^n(t))$ & $(1,0,-1)/n$ & $\{o, \bar{o}, m, \bar{m}\}$ & $1/n$ \\ \hline
\end{tabular}
\end{center}

\end{table}

\emph{The drift, that is, the expected change of $M^n$ in one time step, given the current state of the system is
$f^n(m) = n \mathbb{E} (M^n(t+1)-M^n(t)|M^n(t)=m) $ which can be expressed as:}
$$f^n(m)  =\left(  \begin{array}{c}
-d_{\theta} \delta_D -\! 2 d_{\theta} \delta_m^2 \lambda \frac{n d_{\theta}-1}{n} -\! c_{\theta} \beta \frac{d_{\theta}}{q_{\theta}+d_{\theta}} +\! h(\delta_e \delta_{Sm} {+} \delta_m \eta d) \\
-d_{\theta'} \delta_D -\! 2 d_{\theta'} \delta_m^2 \lambda \frac{n d_{\theta'}-1}{n} -\! c_{\theta'} \beta \frac{d_{\theta'}}{q_{\theta'}+d_{\theta'}} +\! h(\delta_e \delta_{Sm} {+} \delta_m \eta d) \\
2 d_{\theta} \delta_m^2 \lambda \frac{nd_{\theta}-1}{n} - c_{\theta} \delta_C + c_{\theta} \beta \frac{d_{\theta}}{q_{\theta}+d_{\theta}} + h (\delta_{H}+(1-\delta_{H})c  ) \\
2 d_{\theta'} \delta_m^2 \lambda \frac{nd_{\theta'}-1}{n} - c_{\theta'} \delta_C + c_{\theta'} \beta \frac{d_{\theta'}}{q_{\theta'}+d_{\theta'}} + h (\delta_{H}+(1-\delta_{H})c ) \\
d \delta_D {+}  c \delta_C {-} 2h (\delta_{H}+(1-\delta_{H}) c ) {-} 2h(\delta_e \delta_{Sm} {+} \delta_m \eta d    )
 \end{array}  \right)$$
%  \end{array} \!\! \right)}$$
 
\emph{where $m = (d_{\theta},d_{\theta'},c_{\theta},c_{\theta'},h) \,$,  $d =d_{\theta}{+}d_{\theta'} $ and  $c =c_{\theta}{+}c_{\theta'} $. Then the limit of $f^n(m)$ is  }
%is
$$f(m)={\left(\!\!\! \begin{array}{c}
-d_{\theta} \delta_D -\! 2 \lambda d_{\theta}^2 \delta_m^2   -\! c_{\theta} \beta \frac{d_{\theta}}{q_{\theta}+d_{\theta}} +\! h(\delta_e \delta_{Sm} {+} \delta_m \eta d ) \\
-d_{\theta'} \delta_D -\! 2 \lambda d_{\theta'}^2 \delta_m^2   -\! c_{\theta'} \beta \frac{d_{\theta'}}{q_{\theta'}+d_{\theta'}} +\! h(\delta_e \delta_{Sm} {+} \delta_m \eta d   )\\
2 \lambda  d_{\theta}^2 \delta_m^2 - c_{\theta} \delta_C + c_{\theta} \beta \frac{d_{\theta}}{q_{\theta}+d_{\theta}} + h (\delta_{H}+(1-\delta_{H})   c) \\
2 \lambda  d_{\theta'}^2 \delta_m^2 - c_{\theta'} \delta_C + c_{\theta'} \beta \frac{d_{\theta'}}{q_{\theta'}+d_{\theta'}} + h (\delta_{H}+(1-\delta_{H}) c ) \\
d \delta_D {+} c \delta_C {-} 2h (\delta_{H}+(1-\delta_{H})  c   ) {-} 2h(\delta_e \delta_{Sm} {+} \delta_m \eta  d  )
 \end{array}\!\!\! \right)}$$
\emph{Notice that the sum of the all the components of $f(m)$ is zero. Furthermore, if  one of the components  $m_j$ of $m= (d_{\theta},d_{\theta'},c_{\theta},c_{\theta'},h)$  is zero then the corresponding drift function $f_j(m)\geq 0.$ As a consequence, in the absence of birth and death process, the $4-$dimensional simplex is forward invariant, meaning that if  initially $m(0)$ is in the simplex, then for  any time greater than $0$ the trajectory of $m(t)$ stays in the simplex domain.}

\subsubsection{Centralized control design }
\emph{We minimize the proportion of node with states $C$ or $D$  by means of controlling $u(.|),$ i.e., by adjusting $(\delta_m,\delta_e)\in [0,1]^2.$ Since $o(t)=c_1+c_2+d_1+d_2=1-h(t),$ minimizing $o(t)$ is equivalent to maximize the proportion of susceptible node in the population. Therefore the optimization problem becomes}
$$
\left\{
\begin{array}{ll}
\sup_{\ \delta_e,\delta_m}\ & h(T)+\int_0^T h(t) \ dt\\
&\dot{m}=f(m),\ m(0)=m_0\\
& \mbox{where}, \ m=(c_1,c_2,d_1,d_2,h).
\end{array}
\right.
$$

\emph{$\hat{H}= h+ f_1p_1+f_2p_2+f_3 p_3+f_4p_4+f_5p_5.$ This is a twice continuously differentiable function in $m,$
and $\partial_{m_j} \hat{H}=\sum_{i=1}^5 [\partial_{m_j}f_i] p_i $ for $j\leq 4. $ 
The optimum control strategies at time $t$ are the ones that maximize $\hat{H}.$}

$$
\left\{
\begin{array}{l}
\arg\max_{\ \delta_e,\delta_m}\hat{H}\\
\dot{m}=f(m),\ m(0)=m_0\\
\dot{p}_j=-\sum_{i=1}^5 [\partial_{m_j}f_i] p_i,\ j\leq 4,\ t<T\\
\dot{p}_5=-1-\sum_{i=1}^5 [\partial_{h}f_i] p_i , \ t<T\ \\
p(T)=[0,0,0,0,1].
\end{array}
\right.
$$

\subsubsection{Combatting Virus Propagation by Means of Individual Action}
\emph{Let $S(t)$ be the random variable describing the individual state at time $t$ of a generic individual and assume that a generic individual is in a state $s$ at time $t.$ Then $S(t + \frac{1}{n})$ is independent of previous values $(S(t') : t'\leq  t)$ and as $n$ goes  to infinity   for all  state $s'.$
The reward of a generic individual payoff is defined as follows:
$p_{\theta}(s,u,m)=0$ if the individual state $s$ is different than $H,$ and equals $1$ if the state $s=H.$ By doing so, each individual tries to adjust its own trajectory. People in honest state will  accept less meeting and will set their meeting rate 
$\delta_m$ to be minimal, and the other individual with state different than $H$ will try to enter to $H$ as soon as possible. As in a classical communicating Markov chain, this is the entry time to state $H.$}

\emph{Figure \ref{fig:tria} reports the result of the simulation with the following 3 starting points: $(d,c)= (0.2,0.6)$, $(d,c)= (1/3, 1/3)$ and $(d,c)= (0.2,0).$ In the three cases, the system converges to the same steady state  which is around $(d,c)= (0.38,0.6).$ Figure \ref{fig:rew}  plots the reward (honest people) as a function of time for two different control parameters  $\delta_m =0.9$ and $\delta_m =0.1.$ We observe that the reward is greater for $\delta_m =0.1$ than  the one for $\delta_m =0.9.$}

%%%%%

  \begin{figure}[h]
    \centering
    \includegraphics[height=5cm,width=8cm]{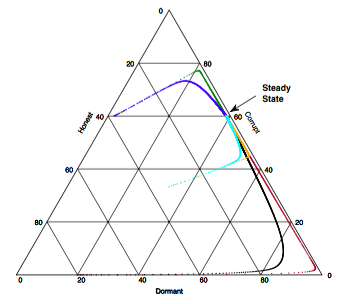} 
 \caption{ {Proportion of dormant, corrupt and honest (followed by the corresponding time-average trajectory).  As time increases, the system approaches a steady state.}
 \label{fig:tria}}
\end{figure}

\begin{figure}[h]
\centering

     \includegraphics[height=5cm,width=10cm]{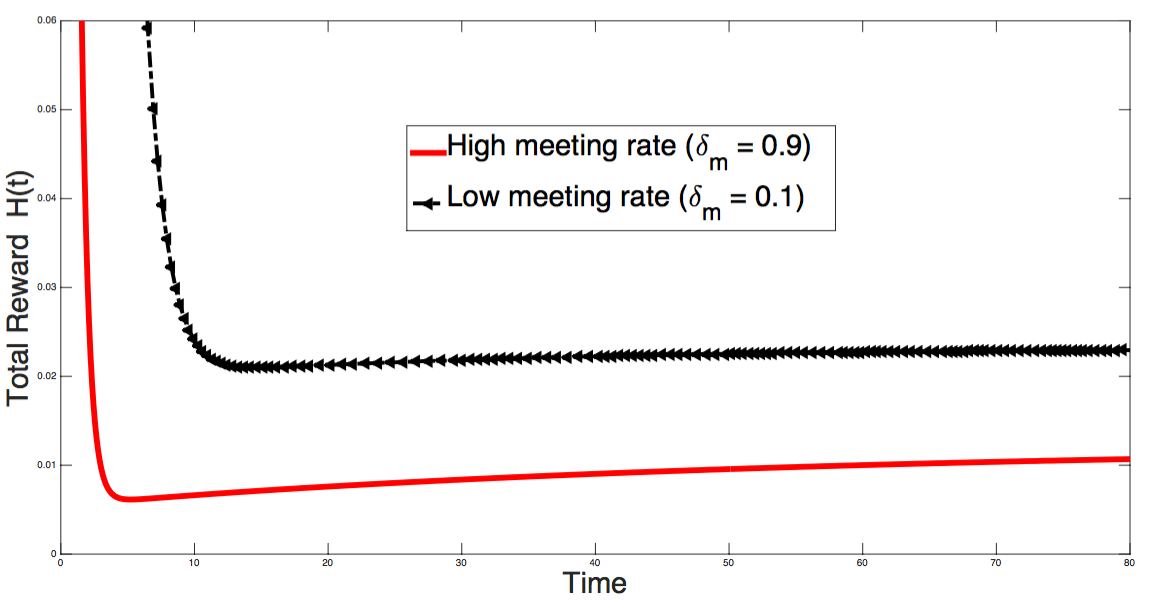}   % eps
\caption{Evolution of Reward (Honest) for the control parameters  $\delta_m =0.9$ and $\delta_m =0.1.$ The smaller the meeting/opening rate is the larger the proportion of  susceptible nodes.}
\label{fig:rew}
\end{figure}

\subsubsection{Network effect}
\emph{The primary advantage of network models is their ability to capture complex individual-level structure in a simple framework. 
To specify all the connections within a network, we can form a matrix from all the interaction strengths which we expect to be sparse with the majority of values being zero. Usually, for simplicity, two individuals (or populations) are either assumed to be connected with a fixed interaction strength or unconnected. In such cases, the network of contacts is specified by a graph matrix $G$, where $G_{ij}$ is 1 if individuals $i$ and $j$ are connected, or 0 otherwise. A connection could be a relationship between the two nodes. It may be represent an internet, social network or physical connection. They may not be close in terms of location. The status of an node will be influenced by the status of its connection following the rules specified above. The resulting graph-based mean-field dynamics is illustrated in Figure \ref{fig:rew3}. }

\begin{figure}[h]
\centering
    \includegraphics[height=5cm,width=9cm]{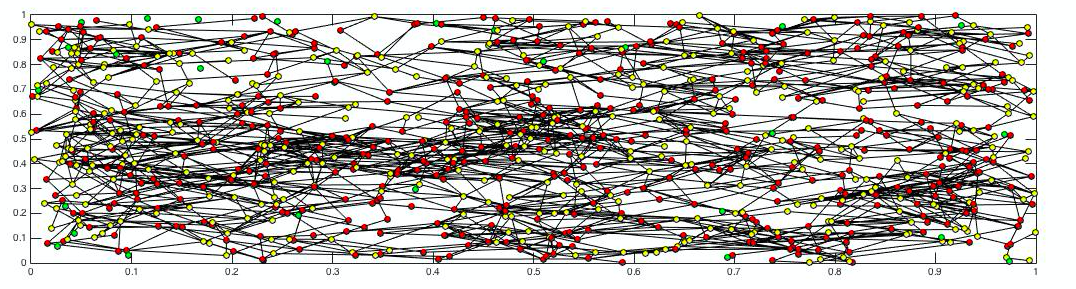} \\
\caption{Network-based virus propagation: each agent has a certain degree of connections without restriction on the location, capturing virus spread via internet or social media contacts. The average degree of the graph is 4. }
\label{fig:rew3}
\end{figure}

\end{appli}

\begin{appli}[Cloud Networks] 
\emph{Resource sharing solutions are very important for data centers as it is required and implemented at different layers of  cloud networks \cite{sha2,cloudnet1,cloudnet2}. The resource sharing problem 
can be formulated as a strategic decision-making problem. Lot of resources may be wasted if the cloud user consider an economic renting. Therefore a careful system design is required when a 
several clients interact. Price design can  significantly improve the resource usage efficiency of large cloud networks. 
 We denote such a game by $\mathcal{G}_{n},$ where $n$ is the number of  clients. The action space of every user is $\mathcal{U}=\mathbb{R}_+$ which is a convex set, i.e.,
each user $i$ chooses an action $u_i$ that belongs to  the set $\mathcal{U}.$ An action may represent a certain demand. All the actions together determine an outcome. Let $p_n$ be the unit price of 
cloud resource usage by the clients. Then, the payoff of user $j$ is given by
\begin{eqnarray}\label{payoffcloud} r_i(x,u_1,\ldots,u_n)=c_n(x)\frac{h(u_i)}{\sum_{j=1}^n h(u_j)}-p_n(x) u_i,  \end{eqnarray} if $\sum_{j=1}^n h(u_j)>0$
and zero otherwise. The structure of the payoff function $r_i(x,u_1,\ldots,u_n)$ for user $i$ shows that it is a percentage of allocated capacity minus the cost for using that capacity.  Here, $c_n(x)$ represents the value of the available resources, $h$ is a positive and nondecreasing function with $h(0)=0.$ We fix the function $h$ to be $x^{\alpha}$ where $\alpha>0$ denotes a certain return index. $x$ is the state of cloud networks which is a random variable on the availability of the servers.
The cloud game $\mathcal{G}_n$ is given by the collection $(\mathcal{X},\mathcal{N}, \mathcal{U}, (r_i)_{i\in\mathcal{I}})$ where $\mathcal{I}=\{1,\ldots, n\},\ n\geq 2,$ is the number of potential users.
The next Proposition provides closed-form expression of the Nash equilibrium of the one-shot game $\mathcal{G}_n$ for a fixed state $x$ such that $c_n(x)>0, p_n(x)>0,$ and for some range of parameter $\alpha.$ It also provides the optimal price $p_n^*$ such that no resource is wasted in equilibrium.}
\begin{proposition} \label{Proposition-NE} \emph{By direct computation,  the following results:}
\begin{itemize}
\item[(1)]\emph{ The resource sharing game $\mathcal{G}_n$ is a symmetric game. All the clients have symmetric strategies in equilibrium whenever it exists.}
\item[(2)] \emph{For  $0\leq \alpha \leq 1,$  and $x\in \mathcal{X},$ the payoff $r_i$ is concave (outside the origin) with  respect to own-action $u_i.$
The best response $BR_i(u_{-i})$ is strictly positive and is given by the root of
$$
z^{(\alpha-1)/2} (\frac{\alpha c_n(x)}{np_n(x)}G)^{1/2}-\frac{z^{\alpha}}{n}-G=0,\ \ G\triangleq\frac{1}{n}\sum_{j\neq i} u_j^{\alpha}
$$
where $z\triangleq u_i$
and there is a unique equilibrium (hence a symmetric one) given by
$ \left(z^{\alpha-1}\frac{\alpha c_n(x)}{n p_n(x)}\frac{n-1}{n} z^{\alpha}\right)^{\frac{1}{2}} -\frac{z^{\alpha}}{n}-\frac{n-1}{n}z^{\alpha}=0,$\ i.e.,
$$u_{NE}^*(x)=\alpha \frac{(n-1)c_n(x)}{n^2p_n(x)}.$$ It follows that the total demand $n a_{NE}^*(x)$ at equilibrium is less than $\frac{c_n(x)}{p_n(x)}$ which means that some resources are wasted. }

\emph{The equilibrium payoff is ${r}_i(x,a_{NE}^*)=u_i p_n(x)\left[ \frac{G+\frac{u^{\alpha}_i}{n}}{\alpha G}-1\right]$ which is positive for $\alpha\leq 1.$}
\item[(3)] \emph{For $\alpha>1,$ the activity (participation) of user $i$ depends mainly of the aggregate of the others. $u_i^*>0$ only if $G\leq G_{*}$ and the number of active clients should be less than $\frac{\alpha}{\alpha-1}.$ If $n>\frac{\alpha}{\alpha-1}$ then $BR_i=0.$}
\item[(4)]\emph{ With a participation constraint, the payoff at equilibrium (whenever it exists) is at least $0.$}
\item[(5)]\emph{ By choosing the price $p_n^*=\alpha\frac{(n-1)}{n}<\alpha$ one gets that the total demand at equilibrium is exactly the available capacity of the cloud. Thus, pricing design can improve  resource sharing efficiency in the cloud. Interestingly, as $n$ grows, the optimal pricing converges to $\alpha.$}
\end{itemize}
\end{proposition}

\emph{We say that the cloud renting game is efficient if no resource is wasted, i.e., the equilibrium demand is exactly $c_n(x).$  Hence, the efficiency ratio is $\frac{n a^*_{NE}}{c_n(x)}.$ As we can see from (ii) of 
Proposition \ref{Proposition-NE}, the efficiency ratio goes to $1$ by setting the price to $p_n^*.$ This type of efficiency loss is due to selfishness and have been widely used in the literature of 
mechanism design and auction theory.  Note that the equilibrium demand increases with $\alpha$, decreases with the charged price and increases with the capacity per user. The equilibrium payoff is positive and  if $\alpha\leq 1$ each user will participate in an equilibrium. In the Nash equilibrium the optimal pricing $p_n^*$ depends on the number of active clients in the cloud and value of  $\alpha.$ When the active number of clients 
varies (for example, due to new entry or exit in the cloud), a new price needs to be setup which is not convenient.}

\end{appli}
\subsection{Mechanical Engineering}

\begin{appli}[Synchronization and Consensus]
\emph{Consider a coupled oscillator dynamics with a control parameter per agent.
$$d\theta_i=[\omega_i+\sum_{j=1}^n K_{ij}(\theta) \sin(\theta_j-\theta_i)  + u_i]dt +\sigma dW_i(t),$$
where $\theta_i$ is the phase of oscillator $i$, $\omega_i$ is the natural frequency of oscillator $i$, $n$ is the total number of oscillators in the system and $K$ is a  coupling interaction term.
The objective here is to explore phase transition and self organization in large population dynamic systems. We explore the mean-field regime of the dynamical mean-field systems and explain how 
consensus and collective motion emerge from  local interactions. These dynamics have interesting applications in multi-robot coordination. Figure \ref{alignment:fig} presents a Kuramoto-based synchronization scheme \cite{KM}. The uncontrolled Kuramoto model can lead to multiple clusters of alignment.  Using mean-field control law, one can drive the trajectories (phases) towards a consensus as illustrated in Figure \ref{fig:consensus} which represents the behaviors for $u_i=-\omega_i+ \eta_i \sin\left( \frac{1}{n}\sum_{j=1}^n \theta_j-\theta_i\right).$ This type of behavior is useful in mobile robot rendezvous problems in which  each agent needs to move towards a common point (where the rendezvous will take place). }

\begin{figure}[h]
\centering

    \includegraphics[height=7cm,width=9cm]{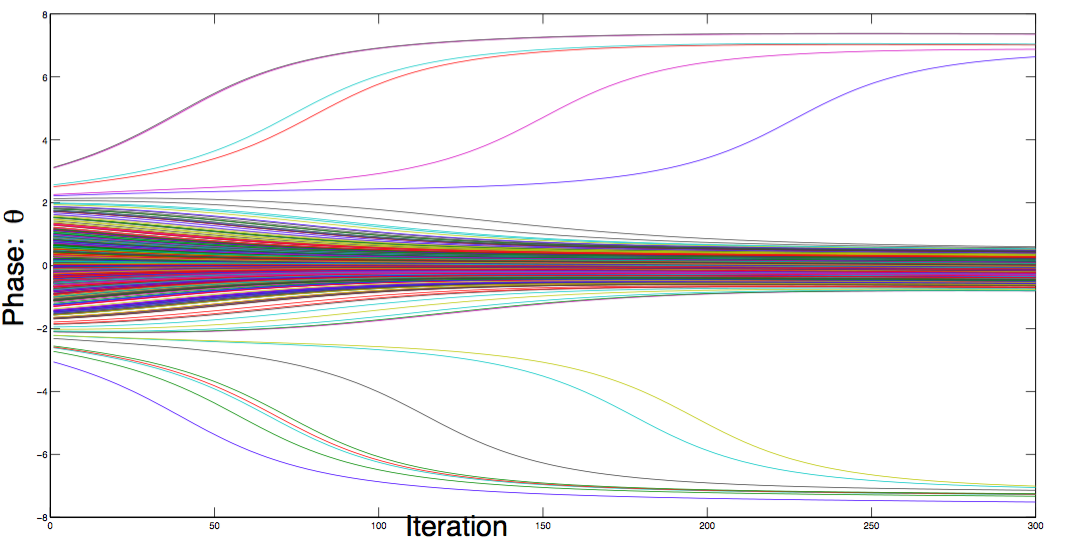} 
\caption{ {Kuramoto-based synchronization scheme with three clusters of alignment with 500 agents.}
\label{alignment:fig}}
\end{figure}

\begin{figure}[h]
\centering

    \includegraphics[height=9cm,width=12cm]{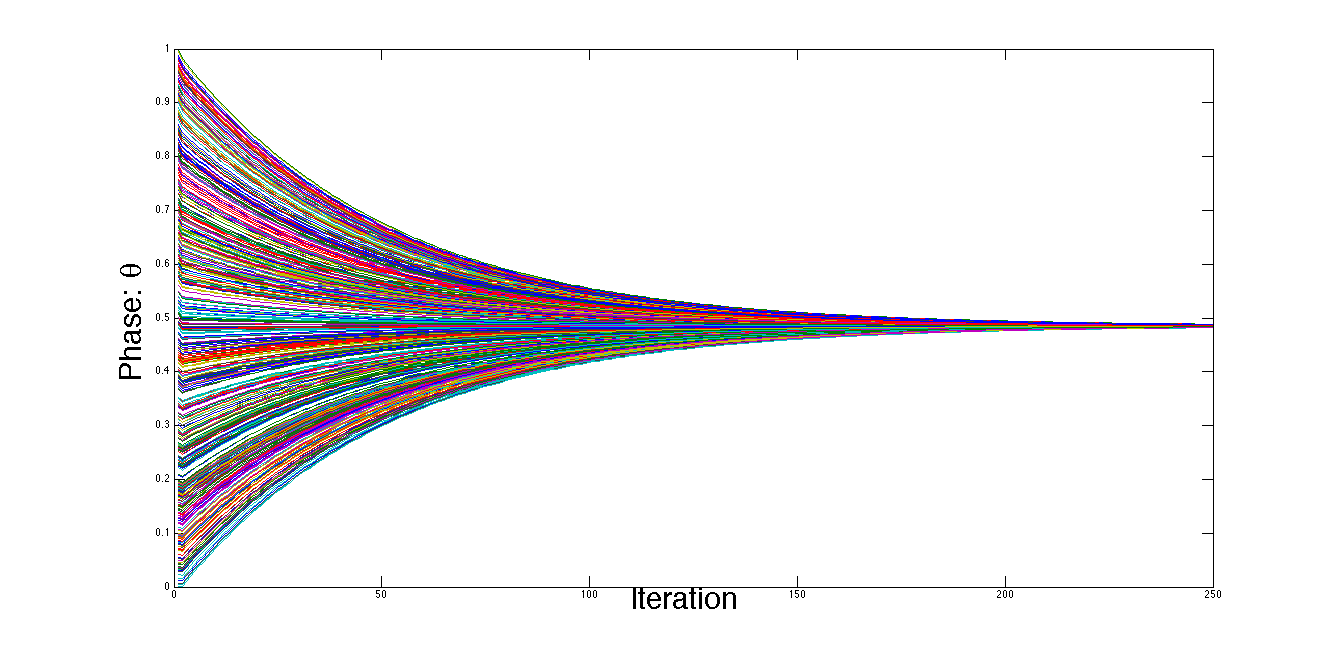} 
\caption{ {A controlled Kuramoto-based synchronization scheme with 500 agents. A mean-field-type control helps to reach a consensus and an agreement independently of the initial distribution of the phases.}
\label{fig:consensus}}
\end{figure}
\emph{ We now provide another relevant application of the Kuramoto model in convoy protection scenario with mobile car-like robots. The goal of the robots is to keep protecting the convoy by occupying the space as the convoy moves. The mean-field-type control helps to balance between energy, placement error and risk.
The authors in \cite{convoy} have shown that  the Kuramoto model   modified with phase shift of $\frac{\pi}{2}$ radians can be used 
in convoy protection scenario given in Figure \ref{fig:auto1}. In this scenario, we want the agents to follow the movement of the convoy while spreading out along a circular perimeter. The mean-field-type control law allow the agents to be positioned equally on a circle and self-organizing the distribution pattern once new agents are added into the network for protecting the convoy and 
occupying the space. Note that re-configuration of the multi-robot team will be done in a distributed way over the circle with center $c$ and with radius $r.$ The protecting convoy is
a rear-wheel drive, front-wheel steerable car-like mobile robot.
 The car-like robot to be controlled is given in Figure \ref{fig:auto2}. The kinematic  parameters of the mobile robot  $i$ are given by 
 $(p_i(t),v_i(t),\theta_i(t),\beta_i(t), l_i)$ representing the cartesian coordinate (position) $p_i(t)=(x_{i,1}(t),x_{i,2}(t))$ of robot $i$ located at the 
 mid-point of the rear-wheel axle, $v_i(t)$ is the translational driving speed, $\theta_i(t)$ is the orientation,
 $\beta_i(t)$
  the steering angle of the front wheels and  $l_i$ the distance between front and rear wheel axle. The goal is to control the robot to a desired orbit while 
  spreading out. One can control the velocity $v_i$ through acceleration and  the steering angle $\beta_i.$ The evolution of center point $c$ and the radius $r$ are given by the drift function 
  $b_c(t), b_r(t).$ The connectivity in the circular graph for agent $i$ is limited to two other agents : $i-1$ and $i+1$ modulo $n.$ Each agent $i$ is influenced only by its neighboring agents.}
\emph{  The instantaneous cost is $$ L_{i}(t)= \epsilon_1[\cos(\theta_{i+1}-\theta_{i})+\cos(\theta_i-\theta_{i-1})]+\epsilon_2(\theta_i-\frac{\pi}{2}-
  \mbox{tan}^{-1}(\frac{x_{i,1}-c_1}{x_{i,2}-c_2})).$$ The first term in bracket says that agent $i$ should spread out from $i-1,i+1.$ The second in the bracket represents  the orientation synchronization.  The terminal cost is of mean-field type and is given by $$ L_{i}(T)=\epsilon_3 \frac{|v_i|^2}{d(x_i,c)^r}+\epsilon_4  |d(x_i,c)-r|^2+\epsilon_5 var(v_i),$$ representing a balance between the kinetic energy spent,  the error adjustment for being on the new circle and the variance respectively.}

\emph{The finite horizon cost functional of agent $i$ is $J_i(u,\beta)=  L_{i}(T)+\int_0^T L_{i}(t) dt.$  Let  $\mathcal{C}(c(0),r(0))$ be the circle with center $c(0)$ and  radius $r(0).$ The best-response problem of agent $i$ is 
 $$ \left\{ \begin{array}{c}
    \sup_{u_i,\beta_i}\ -\mathbb{E}J_i(u,\beta)\\
    \dot{x}_{i,1}=v_i\cos \theta_i,\  \\ 
    \dot{x}_{i,2}=v_i\sin \theta_i,\\
    d{\theta}_i= v_i\frac{tan \beta_i}{l_i}dt+\sigma\frac{tan \beta_i}{l_i} dW_i(t),\\
   v_i= d(x_i,c)^r [ \omega_i+\sum_{j=1}^n K_{ij}(\theta) \cos(\theta_j-\theta_i)  + u_i]\\
       x_i(0)\in \mathcal{C}(c(0),r(0))\subset \mathbb{R}^2
    \end{array} \right.$$
    This is a mean-field-type optimization and the optimality system is easily derived from the stochastic maximum principle. }

\begin{figure}[h]
\centering 
    \includegraphics[height=5cm,width=7cm]{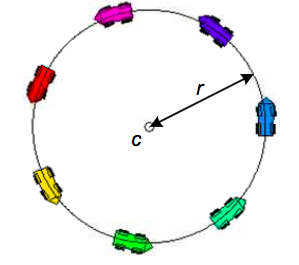} 
\caption{Multi-robot game for protecting a convoy.}
\label{fig:auto1}
\end{figure}

\begin{figure}[h]
\centering
    \includegraphics[height=5cm,width=7cm]{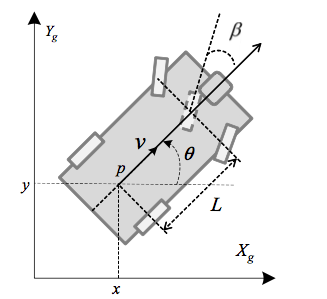} 
\caption{ {Mobile car-like robot.}
\label{fig:auto2}}
\end{figure}
\end{appli}
\begin{appli}[Energy-Efficient Buildings] %
\emph{ Nowadays a large amount of the electricity consumed in buildings is wasted. A major reason for this wastage is inefficiencies in the building technologies, particularly in operating the HVAC (heating, ventilation and air conditioning) systems. These inefficiencies are in turn caused by the manner in which HVAC systems are currently operated. The temperature in each zone is controlled by a local controller, without regards to the effect that other zones may have on it or the effect it may have on others. Substantial improvement may be possible if inter-zone interactions are taken into account in designing control laws for individual zones \cite{pe10, pe11,pe12,pe13,pe14}. The room/zone temperature  evolution is  a controlled stochastic process
$$dT_i= [\epsilon_1(T_{ext}-T_i)+ \sum_{j\in N_i} \epsilon_{2ij}(T_{j}-T_i) + \epsilon_3 u_i (T_{ref}-T_i)]dt+\sigma dW_i,$$ where $\epsilon_1, \epsilon_{2ij}, \epsilon_3$ are positive real numbers. The control action $u_i$ in room $i$ 
depends on the price of electricity { $p(demand, supply,location)$}. The  cost for driving to the comfort temperature zone (see Figure \ref{fig:room})   is  $(T_i-T_{i,comfort})^2+{var(T_i-T_{i,comfort})}.$ The  payoff of consumer  is a sort of 
tradeoff between comfort temperature and electricity cost $u_i p.$ The instantaneous total  cost  of consumer $i$ is} $$L_i(t)= 
\underbrace{u_ip(.)}_\text{energy price}+ \overbrace{(T_i-T_{i,comfort})^2}^\text{deviation to the comfort zone}
+\underbrace{var(T_i-T_{i,comfort})}_\text{risk}.$$

\emph{Within the time horizon $[0,\tau], \ \tau>0,$ consumer $i$ minimizes in }$u_i:$ $$var(T_i(\tau)-T_{i,comfort})+\mathbb{E}\int_0^{\tau} L_i(t) dt.$$ 
\emph{However, the  electricity price $p(.)$ depends on the demand $D=\int_I consumption(i) m_1(t,di)$ and  supply  $S=\int_J supply(j) m_2(t,dj).$ $m_1(t,.)$ is the population mean-field of consumers, i.e., the consumer distribution at time $t.$ Note that $m_1$ is an unnormalized measure.  $m_2$ is the distribution of suppliers.
The building is served by a producer whose remaining energy  dynamics is $$de_{jk}(t)=[c_{jk}(t)\ind_{\{ k \in A^c_j(t)\}} - s_{jk}(t)]dt+\sigma dW_{jk}, $$
 The instant payoff of the producer $j$ is its revenue minus the cost. The cost is decomposed as the cost due to
 mismatch between supply and demand and the production cost. The payoff is  $$r_j=\underbrace{q_j p(D,S)}_\text{revenue}- \overbrace{var(D_j-S_j)}^\text{mismatch cost}- \underbrace{c(q_j)}_\text{production cost}.$$ 
Producer $j$ solves $\max_{q_j} \mathbb{E} \int_0^{\tau} r_j dt$ subject to the production constraint above. 
Explicit solutions to both problem can be obtained using the framework developed in \cite{bou3,andersson}.}

  \begin{figure}[h]
    \centering
    \includegraphics[height=9cm,width=12cm]{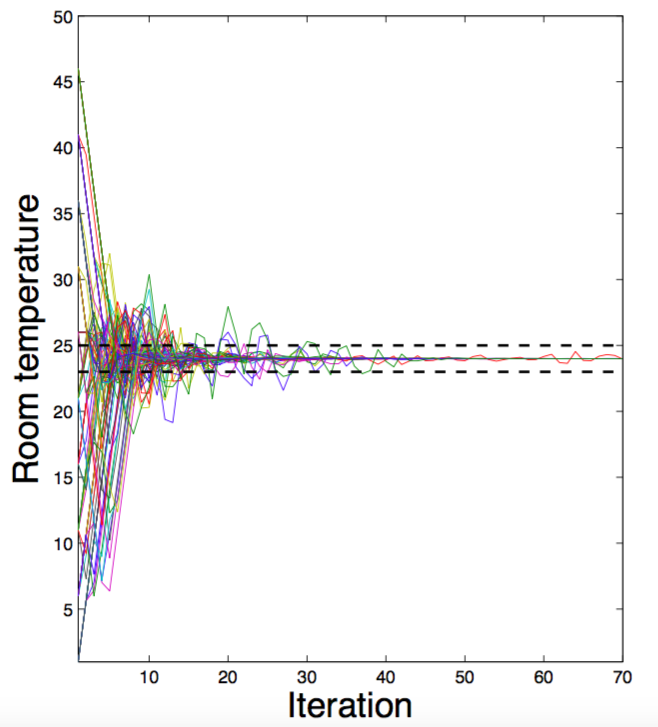} 
 \caption{ {Convergence to comfort temperature between $23$ and $25$ degree celsius  (e.g. 73.4 and 77 Fahrenheit) for 10 connecting rooms in energy-efficient buildings.}}
  \label{fig:room}
  \end{figure}

\end{appli}

\subsection{General Engineering}

 \begin{appli}[ Online Meeting]
\emph{Group meeting online, even over video, is much different than sitting in a boardroom communicating face-to-face with someone. 
But they something in common: deciding to join  Early or on Time the group meeting.
In the context of online video group meeting, since the  communication is over video, the opportunity for miscommunication is much higher, and thus, 
one should pay close attention to how  the group meeting is conducted. Each group member aims to  heighten the quality of her online meetings 
by acting professionally and by signing  early or on time: Nothing throws off a meeting worse than scheduling woes. This is in particular widely observed for online group meetings. }

\begin{figure}[ht]
\centering
 \includegraphics[scale = 0.5]{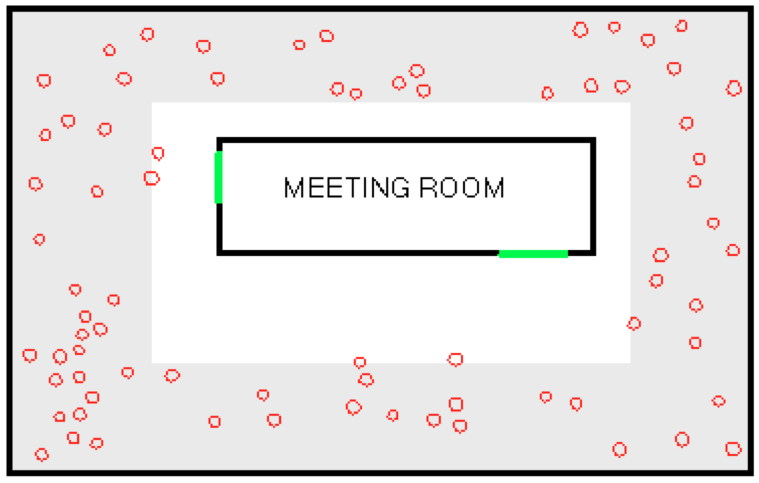}  %
\caption{{Meeting room: initial distribution of the agents represented in 2D} \label{fig:meeting}}
\end{figure}

\emph{Scheduling and synchronization is probably the hardest job in these meetings. The help scheduling groups from different sites can   login to the meeting space at their convenience makes it easier to get meetings started on time.  However, it does not mean that the meeting will start exactly at scheduled time. 
The group members can decide to  be at convenient place early and prepare for the meeting to start, giving you time to settle 
down and get acquainted with the interface. We examine how agents decide when to join the group meeting in a basic setup.
We consider several industry and academia  aiming to collaborate on a research development. The companies are located at different sites. 
Each company from each site has appointed work package leader. In order to improve savings from long business trips, hotels/ accommodation and to reduce jet-lags effect the companies decided to organize an online meeting. After coordinating  all the members availability, date and time is found and the meeting is initially scheduled to start at time $\bar{t}.$ Each member has
 the starting time in his schedule and calendar remainders but in practice, the online meeting only begin when a certain number $\bar{n}$ of representative group leaders and group members  will connect online and will be seated in these respective rooms. Thus, the effective starting time $T$ of the online meeting is unknown and people organize their behavior as a function of 
 $(\bar{t},\bar{n}, T).$ }
 
 {\color{black} \emph{ Each group member can move from her office to the meeting room (see Figure \ref{fig:meeting}). The dynamics of  agent $i$ is simply given by $\dot{x}_i=u_i,$ where $x_i(0)\in D.$ Let  $n(t)$ be the number of people arrived (and seated) in the room before $t.$ If the criterion is met  (by all  groups) before the initially scheduled time $\bar{t}$ of the meeting, this latter starts exactly at   $\bar{t}$. If on the other hand the criterion
   is met at a later time, $T$ is determined by the self-consistency relation: $  T=\inf\{ t\  | \   \ t\geq \bar{t}, \ n(t)\geq \bar{n} \}.$
   The instantaneous cost function is $h(G_n(x_i)) \| u_i\|^2$ and the terminal cost is 
   $c(t_h)=c_1[t_h-\bar{t}]_{+}+c_2[t_h-T]_{+}+c_3[T-{t}_h]_{+}$ where $c_i$ are
    non-negative real numbers, and $ t_h=\inf\{ t,\ | \ x_i(t)\in MeetingRoom  \}.$
    Let $J(u)=c(t_h)+\int_0^{t_h} h(G_n(x_i)) \| u_i\|^2\ dt$  where $h(G_n(x_i)) \| u_i\|^2$  quantifies a congestion-dependent  kinetic energy spent to reach the meeting room of her group. $[T-{t}_h]_{+}$ quantifies  the useless waiting time, $[t_h-T]_{+}$  quantifies of the time for missing of beginning of the online meeting,$[t_h-\bar{t}]_{+}$  quantifies the sensitivity to her reputation of being late at the meeting. Given the strategies  $(u_1,\ldots, u_{i-1},u_{i+1},\ldots, u_n),$ of the other agents, the best response problem of $i$ is:}
  $$ \left\{ \begin{array}{c}
    \sup_{u_i}\ -J(u)\\
    \dot{x}_i=u_i,\  x_i(0)\in D\subset \mathbb{R}^2\\
    u_i=0 \ \mbox{over}\   \partial D\subset \mathbb{R}^2,\
    u_i=k \ \mbox{at}\   \mathrm{Exits} \subset \mathbb{R}^2
    \end{array} \right.$$
    }

\emph{Even if $h(.)$ is constant, the agents interact because of a common term: the starting time of the online meeting $T,$ and 
  $n(T)\geq \bar{n}.$ For this reason, the choice of the other agents matters. 
  The best response of agent $i$ solves the Pontryagin maximum principle }
  $$\begin{array}{c}
  \dot{p}_i= 0, \ t<t_{ih}\\
  \dot{x}_i=u_i^*=\frac{p_i}{2},\\   x_i(0)\in Building \subset  \mathbb{R}^3.
  \end{array}
  $$

\emph{Hence, ${x}_i(t)=x_i(0)+ t  \frac{p_i(t_h)}{2}$ will at arrive at position $x_{room}, $ at time $t_h= 2 \frac{x_{room}-x_i(0)}{p_i(t_h)}$ 
Thus, the optimal payoff of agent $i$ starting from $x$ at time $0$  is  $ -c(t_h)-\int_0^{t_h}  \frac{\|p(t_h)\|^2}{4} dt= -c(t_h)-t_h \frac{\|p(t_h)\|^2}{4}.$
 The optimal payoff of agent $i$ starting from $x$ at time $t$  is  $ -c(t_h)-(t_h-t) \frac{\|p(t_h)\|^2}{4}$
 which is maximized  for  $-c'(t_h)+\frac{\|p(t_h)\|^2}{4}=0,$ i.e., $\|p(t_h)\|^2= 4 c'(t_h)$  hence $\|p(t_h)\|=2\sqrt{c'(t_h)}=\| v_x(t_h, x(t_h))\|.$
  Knowing that the following two functions: $ \tilde{v}_1(x)= \langle x, p^*\rangle,$ with $\| p^*\|_*=1,$  and
$ \tilde{v}_2(x)= c_2\pm \| x-y\|,$ with $x\neq y,$ solves the  Eikonal equation,   $\|\tilde{v}_x\|=1,$  one deduces an
 explicit solution of the Bellman equation:}
$v_t- \frac{\|v_x\|^2}{2}=0,\  \ v(t_h,x)=-c(t_h).$

%can be constructed and it is given by   $${v}(t,x)= -2\sqrt{c'(t_h)} (x(t)-x_{room})- 2(t_h-t)c'(t_h)-c(t_h),$$ and everyone arrives just before the meeting starts. However, the departing time is totally different depending the initially location for the agents. This proves the following result:
\begin{proposition} \label{onlinemeeting:res1}
\emph{The tradeoff value to the meeting room  starting from point $x$ at time $t$ is} ${v}(t,x)= -2\sqrt{c'(t_h)} d(x(t),x_{room})- 2(t_h-t)c'(t_h)-c(t_h).$ 
\end{proposition}

\end{appli}
The next application uses MFTG theoretic modelling for smart cities.
\begin{appli}[Mobile CrowdSensing]
\emph{The  origins of crowdsourcing goes back at least to the nineteenth century and before \cite{cro,cro2}. Joseph Henry, the Smithsonian's first secretary, used the new networked technology of his day, the telegraph, to crowdsource weather reports from across the country, creating the first national weather map of the U.S. in 1856. Henry's successor, Spencer Baird, recruited citizen scientists to collect and ship natural history specimens to Washington, D.C. by the other revolutionary new technology of the day - the railroad - thus forming the bulk of the Institution's early scientific collections.}

\emph{Today's mobile devices and vehicles  not only serve as the key
computing and communication  device of
choice, but it also comes with a rich set of
embedded sensors, such as an accelerometer,
digital compass, gyroscope, GPS, ambient light, dual microphone, proximity sensor, dual camera and many others. %  (see  all the available sensors on iPhone and Samsung Galaxy).  
Collectively, these sensors are
enabling new applications across a wide variety
of domains, creating huge data  and give rise to a new
area of research called mobile crowdsensing or mobile crowdsourcing \cite{cro,cro2,crowdsensing1}. Crowd sensing pertains
to the monitoring of large-scale phenomena that
cannot be easily measured by a single individual.
For example, intelligent transportation systems
may require traffic congestion monitoring and air
pollution level monitoring. These phenomena can
be measured accurately only when many individuals
provide speed and air quality information from
their daily commutes, which are then aggregated spatio-temporally to determine congestion and pollution levels in smart cities.
Such a collected data from the crowd can be seen (up to a certain level) as a knowledge, which in turn, can be seen as a public good \cite{public}. }

\emph{A great opportunity exists to fuse information from populations of privately-held sensors to create useful sensing applications will be public good.  On the other hand, it is important to model, design, analyze and understand the behavior of the users and their concerns such as privacy issues and resource considerations limit access to such data streams. Two MFTGs where each user decides its level of participation to the crowdsensing: (i) public good, (ii) information  sharing, are presented below.}

\emph{The smartphones are battery-operated mobile devices and
sensors suffer from a limited battery lifetime. Hence, there is
a need for solutions that will limit the energy consumptions of
such mobile Internet-connected objects. Such an involvement is translated into a energy consumption cost.}

\emph{ All the data collected from these devices combine both voluntary participator sensing and opportunistic sensing from operators. The data is received by a network of  cloud servers. For security and privacy concerns, several information are filtered, anonymized, aggregated and distributions (or mean-field) are computed.  The model is a public good game with an extra reward for contributors. When decision-makers are optimizing their payoffs, a dilemma arises because individual and social benefits may not  coincide. Since nobody can be excluded from the use of a public good, a user may not have an incentive to contribute to the public good. One way of solving the dilemma is to change the game by adding a second stage in which reward (fair) can be given to the contributors (non-free-riders). }

\emph{The strategic form game with incomplete information denoted by $G_0,$ is described as follows: A stochastic state of the environment is represented by $x.$ There are $n_0$ potential participant to the mobile crowdsensing. The number $n_0$ is arbitrary, and represent the number of users of the game $G_0.$  
As we will see, the important number is not $n_0$ but the number of active users (the ones with non-zero effort), who are contributing to the crowdsensing.}

\emph{ Each mobile user $i$ equipped with sensing capabilities,  can decide to invest a certain level of involvement and effort $u_i\geq 0.$ The action space of user $i$ is $\mathcal{U}_i=\mathbb{R}_+.$ As we will see the degree of participation will be limited so that the action space can be included into a compact interval.
 \color{black} The payoff of user $i$ is additive and has three components: a public good component $\bar G_i(m-\bar R(x)),$ a resource sharing component $\bar R(x) \frac{h_i(u_i)}{\sum_{j=1}^{n_0} h_j(u_j)}$ and a cost component $p(x,u_i).$ Putting together, the function payoff is 
$$  r_{0i}(x,u)=[\bar G_i(m-\bar{R}(x))-p(x,u_i)]\ind_{m\geq \bar R(x)}+\bar R(x) \frac{h_i(u_i)}{\sum_{j=1}^{n_0} h_j(u_j)}\ind_{\sum_{j=1}^{n_0} h_j(u_j)\neq 0}.$$
where $m=\sum_{j=1}^{n_0} u_j$ is the total contribution of all the users, where $\ind_{B}(x)$ is the indicator function which is equal to $1$ if $x$ belongs to the set $B$ and $0$ otherwise. This creates a discontinuous payoff function.  The function $\bar G_i$ is a smooth and nondecreasing, $R(x)$ is a random non-negative number driven by $x.$
The discontinuity of the payoffs due the two branches $\{u:\  m\geq \bar R(x)\}$ and $\{ u: \ m< \bar R(x)\}$ can be handled easily by eliminating the fact that the actions in $\{u: \  m\leq \bar R(x)\}$ cannot be  equilibrium candidates.}

\emph{Using standard concavity assumption with the respect to own-effort, one can guarantee that the game has an equilibrium in pure strategies. We analyze the equilibrium for $\bar G_i(z)=a_i z^{\alpha}, h_i(z)==z$ where $a_i\geq 0,$ and $\alpha\in (0,1].$ 
For any reward $$ \bar{R}(x)\geq \frac{4m^*\sigma}{(1-\sigma)^2} ,\ \sigma=\frac{\bar G'_i(m)-1}{\bar G'_j(m)-1}>0$$ 
where $m^*\in\arg\max [\bar G(m)-m],$ there exists a design parameter $(a_i)_i$
such that the "new" lottery based scheme provides the global
optimum level of contribution in the public good. We collect mobile crowdsensing users to form a network  in which  secondary  users who willing to share their throughput for the benefit of the society or their friends and friends' of friends. 
This can be seen as a virtual Multiple-Inputs-Multiple-Outputs (MIMO) system  with several cells, multiple users per cell, multiple antennas at the transmitters, multiple antennas at the receivers.  The virtual  MIMO system is a sharing network represented by a graph  $(V, E),$ where $V$ is the set of users representing the vertices of the social graph and $E$ is the set of edges. To an active  connection $(i,j)\in E$ is associated a certain value $\epsilon_{ij}\geq 0.$ The term $\epsilon_{ij}$ is strictly positive  if $j$ belongs to the altruistic outgoing network of $i$ and $i$ is concerned about the throughput of user $j.$ The first-order outgoing neighborhood of $i$ (excluding $i$) is $\mathcal{N}_{i,-}.$ Similarly, if $i$ is receiving a certain portion  from $j$ then $i\in \mathcal{N}_{j,-}$ and $\epsilon_{ji}>0.$ In the virtual MIMO system, each user $i$ gets a potential initial throughput $Thp_{i,t}$ during the slot/frame  $t$ and can decide to share/rent some portion of it to its altruism subnetwork members in $\mathcal{N}_{i,-}$. User $i$ makes a sharing decision vector $u_{i,t}=(u_{ij,t})_{j\in \mathcal{N}_i},$ where $u_{ij,t}\geq 0.$ The    ex-post  throughput  is therefore }
$$
Thp_{i,t+}=Thp_{i,t}+\sum_{j \ |\ i \in \mathcal{N}_{j,-}}u_{ji,t}-\sum_{j \in \mathcal{N}_{i,-}}u_{ij,t}.
$$

\emph{Denote $\{ j \ |\ i \in \mathcal{N}_{j,-}\}=: \mathcal{N}_{i,+}.$ Then, }
\begin{eqnarray}
Thp_{i,t+}=Thp_{i,t}+\sum_{j \in \mathcal{N}_{i,+}}u_{ji,t}-\sum_{j \in \mathcal{N}_{i,-}}u_{ij,t}.
\end{eqnarray}
\emph{ Since we are dealing with sharing decisions, the mathematical expressions are not 
necessarily needed if the output can be  observed or measured. Given a measured throughput, A user can decide to share or not based its own needs/demands.
The term $\sum_{j \in \mathcal{N}_{i,+}}u_{ji,t}$ represents the total extra throughput coming to user $i$ from the other users in $N_{i,+}$ (excluding $i$).
The term $\sum_{j \in \mathcal{N}_{i,-}}u_{ij,t}$ represents the total outgoing throughput from user $i$ to the other users in $N_{i,-}$ (excluding $i$). In other  word, user $i$ has shared $\sum_{j \in \mathcal{N}_{i,-}}u_{ij,t}$ to the others. If $j\notin \mathcal{N}_{i,-}$ then $u_{ij,t}=0$ and for all $i,$ $u_{ii,t}=0.$ The balance equation is }
\begin{eqnarray}\nonumber 
\sum_{i}Thp_{i,t+}&=&\sum_{i}Thp_{i,t}+\sum_{i,j}u_{ji,t}-\sum_{i,j}u_{ij,t}\\ &=&\sum_{i}Thp_{i,t},
\end{eqnarray}
\emph{i.e., the system total throughput  ex-post sharing is equal to the system total throughput ex-ante sharing. This means that the virtual MIMO throughput is redistributed and sharing among the users through individual sharing decisions $s.$ Some users may care about the others because he may be in their situation in other slot/day. For these  (altruistic) users, the preferences are better captured by an altruism term in the payoff. We model it through a simple and parameterized altruism payoff. }

\emph{The payoff function of $i$ at time $t$ is represented by }
\begin{eqnarray}
{r}_{1i}(x, u_{i,t}, u_{-i,t})=\hat{r}_i(Thp_{i,t+})+\sum_{j\in \mathcal{N}_i}\epsilon_{ij}\hat{r}_j(Thp_{j,t+}).
\end{eqnarray}

\emph{Here, $\epsilon_{ij}\geq 0$ and represents a certain weight on  how much $i$ is helping $j.$ The matrix $(\epsilon_{ij})$ plays an important role in the sharing game under consideration since it determines the social network and  the altruistic relationship between the users over the network. The throughput $Thp$ depends implicitly  the random variable $x.$
The static simultaneous act one-shot game problem over the network $(V, E)$ is given by the collection
$G_{1,\epsilon}=(V, (\mathbb{R}_+^{n_1-1}, {r}_{1i})_i).$  The vector $u_{i}$ is in $ \mathbb{R}_+^{n_1},$ but the i-th component is $u_{ii}=0.$ Therefore the choice vector  reduces to be in $ \mathbb{R}_+^{n_1-1}.$  and is denoted by $(u_{i,1},\ldots, u_{i,i-1}, 0, u_{i,i+1},\ldots, u_{i,n_1})\ .$ An equilibrium of $G_{1,\epsilon}$ in state $w$ is a matrix $s\in \mathbb{R}_+^{n_1^2}$ such that}
$$
u_{i}\in \mathbb{R}_+^{n_1},\ u_{ii}=0, \ $$ \begin{eqnarray} {r}_{1i}(x,u_{i}, u_{-i})=\max_{u'_{i}} {r}_{1i}(x,u'_{i}, u_{-i}).\end{eqnarray}

\emph{We analyze the equilibria of $G_{1,\epsilon}.$ Note that in practice the shared throughput cannot be arbitrary; it has to be feasible.  }

\emph{The set of actions can be restricted to
$$
\mathcal{U}_i=\left\{ u_i\ | \ u_{ii}=0,\ u_{ij}\geq 0,\ \sum_{j}u_{ij}\leq C\right\},
$$ where $u_i=(u_{i,1},\ldots, u_{i,i-1}, 0, u_{i,i+1},\ldots, u_{i,n}),$ and $C>0$ is large enough. 
For example, $C$ can be taken as the maximum system throughput  $\sum_{j}Thp_{j,0}.$ This way, the set of sharing actions $ \mathcal{U}_i $ of user  $i$ is non-empty, convex and compact. 
Assuming that the  functions $\hat{r}_i$ are strictly concave, non-decreasing and continuous, one obtains  that the game $G_{1,\epsilon}$ has at least one equilibrium  (in pure strategies).}

\emph{As highlighted above, the set of actions can be made convex and compact. Since $\hat{r}_i$ are continuous and strictly convex, it turns out that, each payoff function ${r}_i$  is jointly continuous and is concave in the individual variable $u_i$ (which is a vector) when fixing the other variables. We can apply the well-known fixed-point results which  give the existence of  constrained Nash equilibria.
As we know that $G_{1,\epsilon}$ has  at least one equilibrium, the next step is to characterize them. }

\emph{If  the matrix $u$ is an equilibrium of $G_{1,\epsilon}$ then the following implications hold:}
\begin{eqnarray} u_{ij}>0   \implies \hat{r}'_i(Thp_{i,0+})=\epsilon_{ij}\hat{r}'_j(Thp_{j,0+}).\end{eqnarray}
\emph{The equilibria may not be unique depending on the network topology. This is easily proved and it is due to the fact that one may have multiple ways to redistribute depending on the network structure and several redistributions can  lead to the same sum $Thp_{i,0}+\sum_{j}u_{ji}-\sum_{j}u_{ij}.$ Even if the game has a set of equilibria, the equilibrium throughput and the equilibrium payoff turn out to be uniquely determined. The set of equilibria has a special structure as it is  non-empty, convex and compact. The ex-post equilibrium throughput increases with the ex-ante throughput and stochastically dominates the initial distribution of throughput of the entire network. For $\hat{r}_i=-\frac{1}{\theta}  e^{-\theta Thp_i},\ \theta>0$ let $\epsilon_{ij}=\epsilon$ where $\epsilon>0.$ Then, the  fairness is improved in the network as $\epsilon$ increases. The topology of the network matters.  The difference between the highest throughput  and the lowest throughput in the network is given by  the geodesic distance (strength) of  the multi-hop connection. }
\end{appli}
\section{Time Delayed States and Payoffs}

This section presents MFTGs with time-delayed state dynamics.
Delayed dynamical systems and delayed payoffs appear in many applications. They are characteristic of past-dependence, i.e., their behavior at   time $t$ not only depends on the situation at $t$, but also on their past history and or time delayed state. Some of such situations can be described with controlled stochastic differential delay equations.  Networked systems suffer from intermittent, delayed, and asynchronous communications and sensing. To accommodate such systems, time delays need to be introduced. 

Applications include
\vspace{-3mm}
\begin{itemize}
\item  Consensus and collective motion of Cucker-Smale  \cite{cs} type with delayed information states 
$$
\begin{array}{lll}
dx_i=v_i dt\\
dv_i=\int_{(\bar{x},\bar{v})} a(\| \bar{x}-x_i \|^2)(\bar{v}-v_i)\rho(t-\tau_i, d\bar{x}d\bar{v})\ dt +c \left(\int_{\bar{v}}\bar{v} \rho(t-\tau_i,\mathcal{X}, d\bar{v})\right) \ dt  + u_i dt+\sigma dW_i,
\end{array} $$
where $\rho(t, dxdv)$ is the distribution of states at time $t.$
\item  Delayed information processing, where the difference of the states $\bar{x} -x_i$ influences the dynamics after some time delay $\tau_i$. Examples include Kuramoto-based oscillators \cite{KM}
$$ dx_i=\left[w_i+ \int \rho(t-\tau_i,d\bar{x}) \sin(\bar{x} -x_i(t-\tau_i))+ u_i \right]dt+\sigma dW_i,$$
used to describe synchronization.

\item  Delayed information transmission, where agent $i$ compares its state to the information coming from its neighbor $j $ after some time delay $\tau_i.$ 
Information transmission delays arise naturally in many dynamical processes on networks.  
$$ dx_i=\left[w_i+ \int \rho(t-\tau_i,d\bar{x}) \sin(\bar{x} -{x}_i(t)) + u_i \right]dt+\sigma dW_i.
$$

Delayed information transmission has direct applications in opinion dynamics  and opinion formation on social graph:
$$
dx_i= \left[\int_{B({x}_i,\epsilon_{i})} \rho(t-\tau_i,d\bar{x})  -x_i + u_i \right]dt+\sigma dW_i,
$$
\item The Air Conditioning control towards a comfort temperature is influenced by integrated-state which represents the trend.
\item Transmission and propagation delay affect the performance of both wireline and wireless networks both  delayed information processing and delayed information transmission occur.
\item In computer network security, the proportion of infected nodes  at time $t$ is a function of the delayed state, the topological delay, and the proportion of  susceptible individuals and some time delay for the contamination period. 
\item  In energy markets, there is an observed phenomenon for the dynamics of the price, which comes with a delayed effect.
\end{itemize}

\subsection{Time-delayed mean-field game }
 We consider a  mean-field game where agents interact within the time frame $\mathcal{T}.$ The best-response of a generic agent is 
\begin{eqnarray}
\label{smp0}
\left\{
\begin{array}{lll}
\displaystyle{\sup_{u\in \mathcal{U}} \mathbb{E} \left[ G(u,m_1,m_2)\right] },
\
\displaystyle{\mbox{ subject to }\ }\\
dx= b(t, x, y,z,u,m_1,m_2,\omega) dt \\ + \sigma(t, x, y,z,u,m_1,m_2,\omega) d{W}\\
+\int_{\Theta}\ {\gamma}(t,x,y,z,u,m_1,m_2,\theta,\omega) \tilde{N}(dt,d\theta),\\
x(t)=x_{0}(t),\ \   \ t\in [-\tau,0], \\
 \end{array}
\right.
\end{eqnarray}
where  $\tau_k>0$ represents a time delay, $x=x(t)$ is the state at time $t$ of a generic agent,
$ y=(x(t-\tau_k))_{1\leq k\leq D}, $  is a $D-$dimensional delayed state vector,
 $ z(t)=(\int_{t-\tau}^t \lambda(ds) \phi_l(t,s)x(s))_{l\leq I}$ is the integral state vector of the recent past  state over $[t-\tau,t].$  This represents the trend of the state trajectory. The process $\phi_l(t,s)$ is an $\mathcal{F}_s-$adapted locally bounded process.  $\lambda$ is a positive and $\sigma-$finite measure. $m_1$ the average states of all the agents, $m_2$ the average control actions of all the agents, $x_0$ is a initial deterministic function of  state.  $W(t)=W(t,\omega)$ be a standard Brownian motion  on $\mathcal{T}=[0,T]$ defined on a given filtered probability space $ (\Omega, \mathcal{F}, \mathbb{P}, \{\mathcal{F}_t\}_{t\in \mathcal{T}}).$

%%%%%%%%%%%%%%%%%%%%%%%%%%%%%
%%%%%%%%%%%%%%%%%%%%%%%%%%%%%%%%%%%%%%%%%%%%%
 {Payoffs: }  $G(u,m_1,m_2)= g_1(T,x(T), m_1(T),\omega){+}\int_{t\in \mathcal{T}} g_0(t, x, y,z,u,m_1,m_2,\omega)\ dt,$ where the instantaneous payoff function is $g_0:\ \mathcal{T}\times \mathcal{X}^3\times 
{U}\times \mathcal{X}\times U\times \Omega\rightarrow \mathbb{R},$ the terminal payoff function is $g_1:\ \mathcal{X}^2 \times \Omega\rightarrow
\mathbb{R}.$  

 { State dynamics:}
The drift  coefficient function is $b: \ \mathcal{T}\times \mathcal{X}^3\times 
{U}\times \mathcal{X}\times U\times \Omega \rightarrow \mathbb{R},$   the  diffusion coefficient function is  $ \sigma:\ \ \mathcal{T}\times \mathcal{X}^3\times 
{U}\times \mathcal{X}\times U\times \Omega \rightarrow \mathbb{R}.$ 
 
 {Jump process:}
Let ${N}$ be a Poisson
random measure with L\'evy measure $\mu (d\theta),$  independent of $\mathcal{B}$ and the measure $\mu$ is a $\sigma-$finite measure over $\Theta.$ 
$\tilde{N}(dt,d\theta)=N(dt,d\theta)-\mu(d\theta)dt.$
The function ${\gamma}: \ 
 \mathcal{T}\times \mathcal{X}^3\times 
{U}\times \mathcal{X}\times U\times \Theta\times \Omega   \rightarrow \ \mathbb{R}. $ The filtration $\mathcal{F}_t$ is the one generated by the union of events from ${W}$ or ${N}$ up time $t.$

The goal is to find or to characterize  a best response strategy to mean-field $(m_1,m_2): $
$u^*\in  \arg\max_{u\in\mathcal{U}} G(u, m_1, m_2).$

{\bf Hypothesis H1: } The functions $b, \sigma, g$ are   continuously differentiable with the respect to $(x,m).$ Moreover, $b, \sigma, g$ and all their first  derivatives with the respect to $(x,y,z,m)$ are 
continuous in $(x,m, u)$ and bounded.

We explain below why the existing solution approaches cannot be used to solve  $(\ref{smp0}).$ First, the presence of $y,z$ lead to a delayed integro-McKean-Vlasov  and the stochastic maximum principle developed  in \cite{alainv1,PLref3,peter,carmona,pri1,pri4} does not apply. The dynamic programming principle for Markovian mean-field control cannot be directly used  here because the state dynamics is non-Markovian  due to the past  and time delayed states. 
Hence, a novel solution approach or an extension is needed in order to solve $(\ref{smp0}).$  A chaos expansion methodology can be developed as in  \cite{cdc15} using generalized polynomial of Wick and  
Poisson jump process. The idea  is to develop a finite-dimensional optimality equation for $(\ref{smp0}).$ In this respect, a stochastic maximum principle could be a good candidate solution approach.
Under H1, for each control $u\in \mathcal{U},$  $m_1$ and $m_2$ the state dynamics admits a unique solution, $x(t):=x^{u}(t).$
The non-optimized Hamiltonian  is
$H(t,x,y,z,u,m_1,m_2, p,q,\bar{r},\omega):   \mathcal{T}\times \mathcal{X}^3\times 
{U}\times \mathcal{X}\times U\times \mathbb{R}^2\times J \times \Omega \rightarrow \mathbb{R}$
where $\bar{r}(.)\in J$ and $J$ is the set of functions on $\Theta$ such that $\int_{\Theta} \gamma \bar{r}(t,\theta)\mu(t,d\theta)$  is finite.
The Hamiltonian is $H=g_0+b p+\sigma q+\int_{\Theta} \gamma \bar{r}(t,\theta)\mu(d\theta).$
The first-order adjoint process  $(p,q,\bar{r})$ is  time-advanced and determined by
\begin{eqnarray}\nonumber dp&=&E[-H_x \ind_{ t\leq T} -\sum_{k=1}^D H_{y_k}(t+\tau_k)\ind_{ t\leq T-\tau_k} \ | \ \mathcal{F}_t] dt\\ \nonumber
&& -\sum_{l=1}^I E[\lambda(dt)\int^{t+\tau}_t\phi_l(t,s)  H_z \ind_{s\in [0,T]} ds\   \ | \ \mathcal{F}_t]\\ \label{foad2}
&& +qdW(t)
+\int  \bar{r}(t,d\theta) \tilde{N}(dt,d\theta),\\ \label{foad}
p(T)&=&g_{1,x}(x(T), m_1(T)).\  \ 
\end{eqnarray}

We now discuss the existence and uniqueness of the  first-order adjoint equation. 
{\color{black}
\begin{proposition} \label{mainr1} \emph{Assume that the coefficients are $L^2,$ the first order adjoint (\ref{foad}) has a unique solution such that
$$  \mathbb{E}\left[   \int_0^T p^2+q^2+\int_{\Theta} \bar{r}^2(t,\theta)\mu(d\theta) \ dt\right] <+\infty$$
Moreover, the solution $(p,q,\bar{r})$ can be found backwardly as follows:}
\begin{itemize}
\item \emph{Within the time frame $(T-\tau, T),$    $dp= E[-H_x \ | \ \mathcal{F}_t] dt+qdW(t)
+\int_{\Theta}   r(t,d\theta) \tilde{N}(dt,d\theta)$ with $p(T).$}
\item \emph{We fix $p(T-\tau)$ from the previous step and solve  (\ref{foad2}) on interval $(T-2\tau, T-\tau).$}
\item \emph{We inductively construct a procedure to compute  $p(t)$ on $t\in [T-k\tau, T-(k-1)\tau],\ k\leq \frac{T}{\tau}$ ending with  $p(T-(k-1)\tau).$   }
\end{itemize}
\end{proposition}
}

Note that,  if $t\in [T-k\tau, T-(k-1)\tau]$  then  $t+\tau\in [T-(k-1)\tau, T-(k-2)\tau]$ and  hence, $(p(t+\tau),q(t+\tau),\bar{r}(t+\tau,\theta))$ is known from the previous step. However, $p(t+\tau)$ may not be $\mathcal{F}_t-$adapted. Therefore a conditional expectation with the respect to the filtration $\mathcal{F}_t$ is used.

If $U$ is a convex domain, we know that the second-order adjoint processes of Peng's type are not required, and if $(x^*, u^*)$ is a best response 
to $m_1,m_2$ then there is a triplet of processes $(p,q,\bar{r}),$ that satisfy the first order  adjoint equation such that
\begin{eqnarray}H(t, x^*,y^*,z^*, u^*, m_1,m_2, p,q,\bar{r})  \nonumber \\ -H(t, x^*,y^*,z^*, u, m_1,m_2, p,q,\bar{r})
  \geq 0,\end{eqnarray} for all $u\in \mathcal{U},$ almost every $t$ and $\mathbb{P}-$almost surely (a.s.).
A necessary condition for (interior) best response strategy is therefore $E[ H_u\ \ | \ \mathcal{F}_t]=0$ whenever $H_u$ makes sense. A sufficient condition for optimality can be obtained, for example, in the concave case: $g_1, H$ are concave in $(x,y,z,u)$ for each  $t$ almost surely.
  
  \subsection{Time delays effect in the Prosumers' Integration to Power Networks}
  
{\color{black} Let  $c_1(t), c_2(t)$ and $c_3(t, z) $ be given bounded adapted processes, with $c_1$ assumed to be deterministic and $\int c_3^2 \nu(dz) <+\infty.$
Consider the energy dynamic generated by a prosumer  as
$$
de_{i}=   (c_1(t) e_i(t-\tau) -u_i)dt +c_2(t)e_i(t-\tau) dW(t) +e_i(t-\tau)\int c_3(t,\theta)\tilde{N}(dt,d\theta),
$$

$e_i(t)=e_{i0}(t)\ind_{[-\tau, 0]}(t)$ where $e_{i0}$ is deterministic and bounded function that is given. The  energy $u_i$ is consumed by $i.$ Prosumer $i$ has a (random) 
satisfaction function $s(t,u_i,\omega)$ which is $\sigma(W_{t'}, N(t'),\ t'\leq t)-$adapted for each consumption strategy $u_i\geq 0,$ the random function $s$ is assumed to be continuously differentiable and increasing with the respect to $u_i$ and its derivative $s_{u_i}(t,u_i,\omega)$ is decreasing in $u_i$. The function $s_{u_i}(t,u_i,\omega)$ vanishes as the consumption  $u_i$ grows without bound. Therefore, the maximum value of $s_{u_i}(t,u_i,\omega)$ is achieved when $u_i=0$ and the maximum value is 
$\bar{m}(t,\omega):=s_{u_i}(t,0,\omega).$ The infinimum value of $s_{u_i}(t,u_i,\omega)$ is $0.$
   It follows that  $u_i \mapsto s_{u_i}(t,u_i,\omega)$ is a one-to-one mapping from $\mathbb{R}_{+}$ to  $(0, \bar{m}(t,\omega)].$ In particular, the function $br: \ \lambda \mapsto (s_{u_i}(t,.,\omega))^{-1} [\lambda]\ind_{(0, \bar{m}(t,\omega)]}(\lambda)$  is well-defined and is a measurable function.
}
Prosumer $i$ aims to maximize her satisfaction functional together with her profit $\mathbb{E}\left[ g(e_i(T))+ \int_0^T s(t,u_i,\omega)+price(m) q_i\ dt \right]$ 

The Hamiltonian is 

$H(t, x,y,z, u_i, m_1,m_2, p,q,\bar{r})=s+(c_1 y -u_i) p+c_2y q+ y\int_{\Theta} c_3 \bar{r}(t,\theta)\mu(d\theta).$

\begin{eqnarray}\nonumber dp&=&E[- H_{y}(t+\tau)\ind_{ t\leq T-\tau} \ | \ \mathcal{F}_t] dt +qdW(t)
+\int  \bar{r}(t,d\theta) \tilde{N}(dt,d\theta),\\ \label{foad3}
p(T)&=&g_{x}(x(T)),\  \ 
\end{eqnarray}
where 
$H_{y}(t+\tau)= c_1(t+\tau)p(t+\tau)+c_2(t+\tau)q(t+\tau)+ \int_{\Theta} c_3(t+\tau) \bar{r}(t+\tau,\theta)\mu(d\theta).
$

We solve the solution explicitly with $g(x)= c_4 x,\  \ c_4\geq 0.$
$p(T)=c_4\geq 0.$ Between time $T-\tau$ and $T,$ the stochastic process $p(t)$ must solve
$dp=qdW(t)
+\int  \bar{r}(t,d\theta) \tilde{N}(dt,d\theta)$ and it should be $ \mathcal{F}_t$-measurable. Therefore $p(t)=c_4$ on $t\in [T-\tau, T].$
For $t< T-\tau,$ the processes $q $ and $\bar{r}$ are zero and $p$ is entirely deterministic and solves
$$
\dot{p}=-c_1(t+\tau)p(t+\tau).
$$
Thus, for $t\in [T-2\tau, T-\tau],$ 
$$
p(t)=p(T-\tau)+\int^{T-\tau}_t c_1(t'+\tau)p(t'+\tau) \ dt'.$$
This means that $p(t)= c_4[1+\int^{T}_{t+\tau} c_1(t'') \ dt''].
$ For $t\in [T-(k+1)\tau,T-k\tau],$  and $(k+1)\tau\leq T,$ one has
$
p(t)=p(T-k\tau)+\int^{T-(k-1)\tau}_{t+\tau} c_1(t'')p(t'') \ dt''.
$

By assumption, $s_{u_i}(t,u_i,\omega)$ is decreasing in $u_i$ and from the above relationship it is clear that $p$ is decreasing with $\tau.$ It follows that, if $\tau_1<\tau_2,$
$p[\tau_1](t) > p[\tau_2](t).$ We would like to solve   $s_{u_i}(t,u_i,\omega)=p[\tau_1](t) > p[\tau_2](t).$
By inverting the above equation one gets
$u_i^*[\tau_1] < u_i^*[\tau_2].$
Thus, the optimal strategy $u_i^*$ increases if the time delay $\tau$ increases.

This proves the following result:
\vspace{-2mm}
\begin{proposition} \label{delay:res2}
\emph{The time delay decreases the prosumer  market price.  The optimal strategy $u_i^*$ increases as the time delay $\tau$ increases. }
\end{proposition}

\vspace{-2mm}
Numerical methods for delayed stochastic differential equations of mean-field type is not without challenge.  Here we implement the  Milstein scheme using MATLAB. We choose  the following parameters $\gamma=0,c_1=c_2=c_3=1$ and  set  the satisfaction function as $$s(u)=1-(1+\mu\bar{m}_2)e^{-u}$$ where $\mu>0$ and $\bar{m}_2$ is the average of all other agents' control actions. A typical shape of the satisfaction function is given in Figure \ref{shape1}. The optimal control is $$u^*(t)=- \frac{\log p(t)}{1+\mu \bar{m}_2(t)}\ind_{(0,1]}(p(t)).$$

$$
u^*(t)=\left\{\begin{array}{lll}
-\frac{\log c_4}{1+\mu \bar{m}_2(t)} \ \mbox{on}\ t\in (T-\tau, T],\\  \\ 
-\frac{\log c_4(1+T-t-\tau)}{1+\mu \bar{m}_2(t)}  \  \mbox{on}\ t\in (T-2\tau, T-\tau],\\  \\ 
-\frac{1}{1+\mu \bar{m}_2(t)}\log [c_4(1+\tau) +c_4(1+T-\tau)(T-t-2\tau)  -\frac{c_4}{2}(T-t-2\tau)(T+t)]  \\ \mbox{on}\ t\in (T-3\tau, T-2\tau].
\end{array} 
\right.
$$

The mean-field equilibrium solves the fixed-point equation $\mathbb{E}[u^*(t)]=\bar{m}_2(t).$ Putting together one obtains
\vspace{3mm}
$$\bar{m}_2(t)=- \frac{\log p(t)}{1+\mu \bar{m}_2(t)},$$ i.e., the root (in $\bar{m}_2$)  of $\bar{m}_2 \mapsto \bar{m}_2(1+\mu \bar{m}_2)+\log p(t).$ The quadratic polynomial has two roots: one positive and the other negative value. Since the consumption is nonnegative, the mean of the mean-field control action  is hence given by 
$$
\bar{m}_2(t)=\frac{-1+\sqrt{1+4\mu\log [\frac{1}{p(t)}]}}{2\mu}.
$$

Notice that the effect of the time delay  $\tau$ in this specific example was through the adjoint process $p$ which also enters into the control action $u$.
\begin{figure}[htb]
\centering
\includegraphics[width=0.5\textwidth]{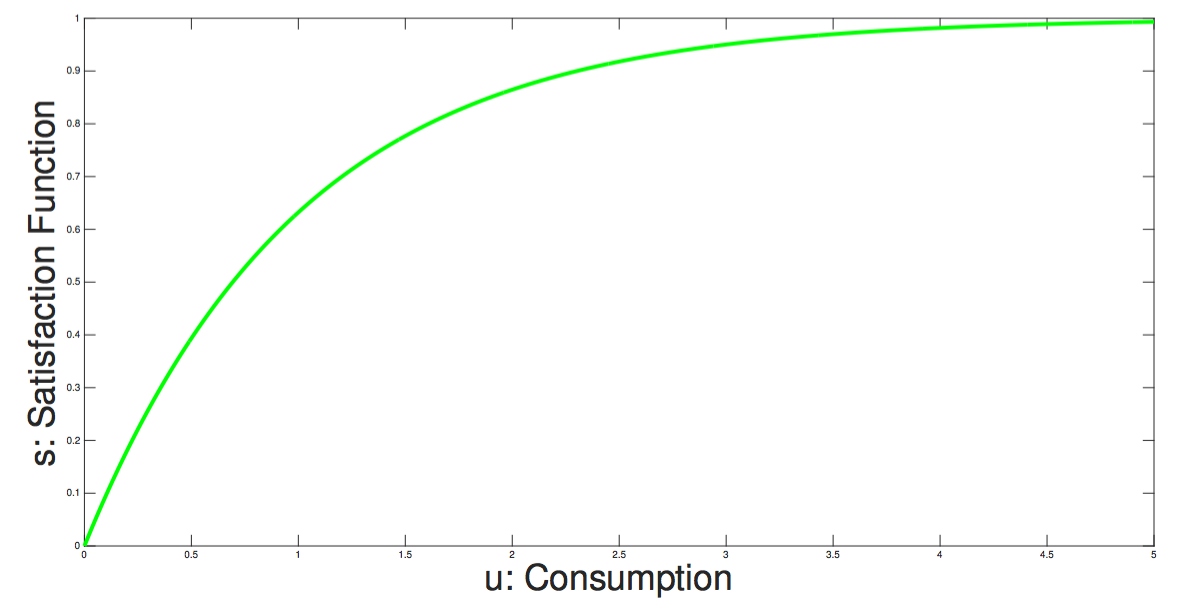}\\
\caption{Typical shape of the satisfaction function of the prosumer.}
\label{shape1}
\end{figure}

\begin{figure}[htb]
\centering
\includegraphics[width=0.5\textwidth]{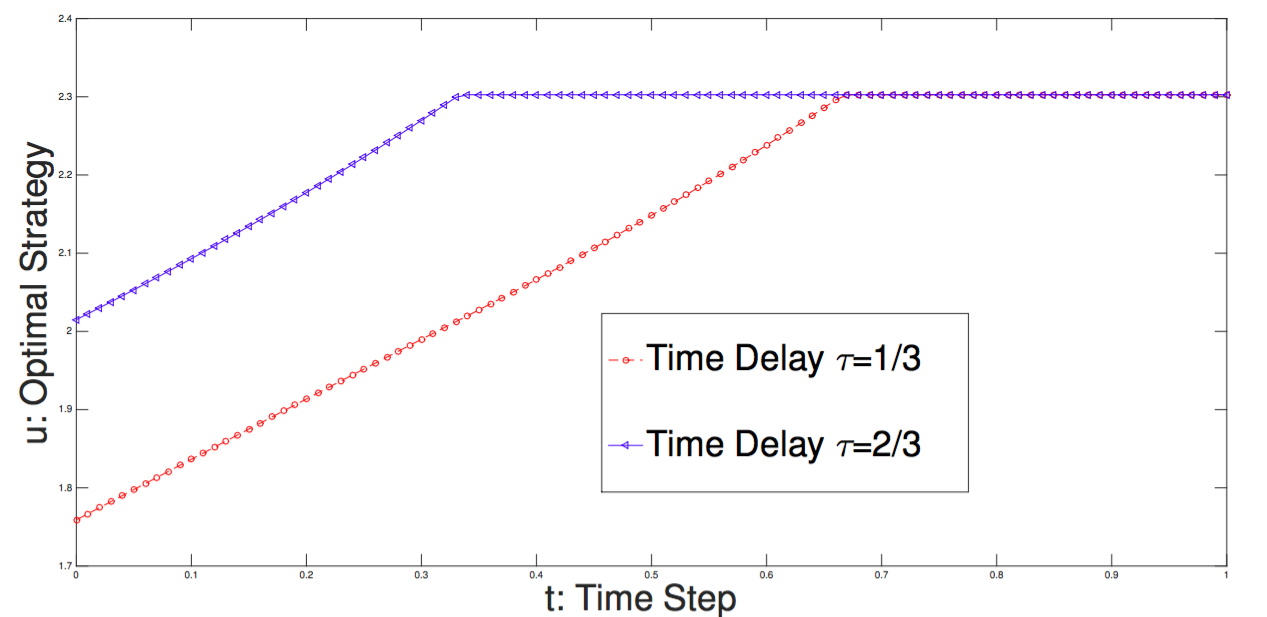}\\
\caption{Structure of the Optimal Strategy for $T=1,\tau=1/3, \tau=2/3.$}
\label{optcontrol}
\end{figure}\
\begin{figure}[htb]
\centering
\includegraphics[width=0.5\textwidth]{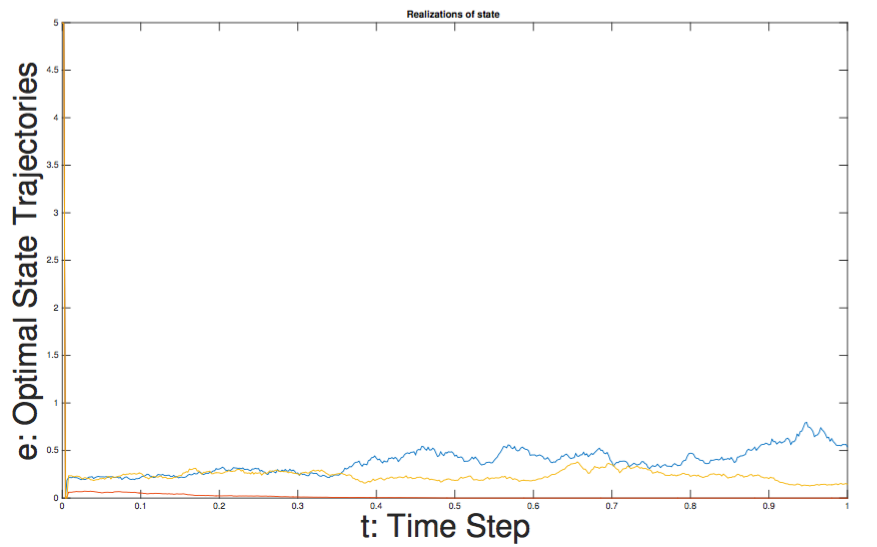}\\
\caption{Sample optimal state trajectories for $T=1,\tau=1/3$ using Milstein scheme.}
\label{sample}
\end{figure}\
We plot the structure of the optimal strategy for $T=1,\tau=1/3, \tau=2/3.$ The theoretical result of Proposition 2 is numerically observed in Figure \ref{optcontrol}. Figure \ref{sample} plots sample optimal state trajectories for $T=1,\tau=1/3$ using Milstein scheme. 

\section{Decentralized Information  and Partial Observation}
{\color{black}Let  $\mathcal{F}_t^{W}$  be the $\mathbb{P}$-completed natural filtrations generated by $W$ up to $t$. 
Set $\mathcal{F}^W:=\{\mathcal{F}_t^{W},\ 0\leq t \leq T\}$ and $\mathbb{F}:=\{{\mathcal{F}}_t,\ 0\leq t \leq T\}$, where $\mathcal{F}_t=\mathcal{F}_t^{W} \vee \sigma(x_0)$. 
An admissible control $u_i$ of agent $i$  is an $\mathcal{F}^{W_i}$-adapted process with values in a non-empty, closed and bounded subset (not necessarily convex) $U_i$ of $\mathbb{R}^{d}$ and satisfies $E[\int_0^T|u_i(t)|^2dt]<\infty$. Those are nonanticipative measurable functionals of the  Brownian motions.
Since each agent  has a different information structure (decentralized information), let
$\mathcal{U}_i$ be the set of admissible strategies of $i$ (with decentralized partial information) such that  $\mathcal{G}_{i,t} \subset \mathcal{F}_{i,t},$ i.e., 
$\mathcal{U}_i:= \{ u_i\in  L^2_{\mathcal{G}_{i,T}}([0,T], \mathbb{R}^{d}),\ u_i(t,.)\in U_i \  \mathbb{P}-a.s\}.$
Given a strategy $u_i\in \mathcal{U}_i$, and a (population) mean-field term $m$ generated by other agents we consider the signal-observation $x_i^{u_i,m}$ which satisfies the following stochastic differential equation of mean-field type to which we associate a best-response to mean-field \cite{bou1,bou2,bou3}:
  \begin{eqnarray}
\label{type2}   
\left\{ 
\begin{array}{lll} \ \sup_{u_i\in  \mathcal{U}_i} R(u_i,m)\
\displaystyle{\mbox{ subject to }\ }\\
dx_i(t)= b(t, x_i(t), Ex_i(t), u_i(t), m(t)) dt +\sigma(t, x_i(t), Ex_i(t), u_i(t), m(t)) dW_{i,t} ,\\
x_i(0)\sim \mathcal{L}(X_{i,0}),\\
m(t)=\mbox{population mean-field },
 \end{array}
\right.
\end{eqnarray}

\begin{equation}
b(t,x,y,u,m): \,\,[0,T] \times \mathbb{R}^d\times\mathbb{R}^d\times U_i \times \Lambda \longrightarrow \mathbb{R},
\end{equation}
\begin{equation}
\sigma(t,x_i,y_i,u_i,m): \,\,[0,T] \times \mathbb{R}\times \mathbb{R}\times U_i \times \Lambda\longrightarrow \mathbb{R}.
\end{equation}
$$
\begin{array}{l} \
R(u_i,m)= g(x_i(T), Ex_i(T), m(T)) +\int_0^T r(t, x_i(t), Ex_i(t), u_i(t), m(t)) dt,
\end{array}$$  $g$ is the terminal payoff and $r$ is the running payoff.}
Given $m,$ any $ u^*_i\in {\cal U}_i$ which satisfies $R( u_i^*(\cdot),m)=\sup_{u_i(\cdot)\in {\cal U}_i}R(u_i, m)$
is called a {\it pure} best-response strategy to $m,$ by agent $i.$
In addition to the other coefficient we assume that $\gamma$ satisfies H1. Under H1, the state dynamics  admits a unique strong solution (see \cite{jourdain}, Proposition 1.2.)
Given $m,$ we apply the  SMP for risk-neutral mean-field type control from  (\cite{b3}, Theorem 2.1) to the state dynamics $x$ to derive the first order adjoint equation.
Under the assumption H1, there exists a unique $\mathbb{F}$-adapted pair of processes $(p, q)$, which solves the Backward SDE:
\begin{eqnarray}\label{rn-firstAD} p(t)&=&g_x(T)+E[g_y(T)]\\  && \nonumber +\int_t^T \{ H_x(s)+E[H_y(s)]\}ds\\  && \nonumber -\int_t^T q(s) dW(s),
\end{eqnarray}
such that  $ \mathbb{E}\left[\sup_{t\in [0,T]}\ |p(t)|^2+\int_0^T |q(t)|^2dt  \right] < +\infty.$
However, these processes $(p,q)$ may not be adapted to decentralized information  $\mathcal{G}_{i,t}.$ This is why their conditioning will  appear in the maximum principle below.
Again by (\cite{b3}, Theorem 2.1), there exists a unique $\mathbb{F}$-adapted pair of processes $(P, Q)$, which solves the second order adjoint equation 
\begin{eqnarray}\label{rn-secondAD}  && P(t)= g_{xx}(T)\\  && \nonumber + \int_t^T \{2 b_x(s) P(s)+\sigma_x^2 P(s) +2 \sigma_x(s) Q(s)+H_{xx}(s)\} ds \nonumber \\ 
&&-\int_t^T Q(s) dW(s),
\end{eqnarray}
such that  $ \mathbb{E}\left[\sup_{t\in [0,T]}\ |P(t)|^2+\int_0^T |Q(t)|^2dt  \right] < +\infty.$
Note that in the multi-dimensional setting, the term $2 b_x(s) P(s)+\sigma_x^2 P(s) +2 \sigma_x(s) Q(s)$ becomes
 $b'_x P+  Pb_x+\sigma'_x P \sigma_x + \sigma'_x Q+Q\sigma_x.$
 \begin{proposition} \label{rn-SMP2}   \label{mainr2}
\emph{Let H1 holds and $m$ be a given population mean-field profile. If $( x^*_i, u^*_i)$ is a best-response then, there are two pairs of $\mathbb{F}$-adapted processes $(p,q)$ and $(P, Q)$ that satisfy (\ref{rn-firstAD}) and (\ref{rn-secondAD}) respectively, such that}
\begin{equation}\label{rn-VI2}\begin{array}{lll}
i\in \mathcal{N}: \ \left[\delta H(t)+ \frac{1}{2}\delta\sigma(t)'P(t)\delta\sigma(t) \  | \ \mathcal{G}_{i,t}\right] \leq 0, 
\end{array}
\end{equation} \emph{for all $u_i\in  \mathcal{U}_i,$ almost every $t$ and $\mathbb{P}-$almost surely,
where,}
\begin{equation}\label{rn-VItt2}\begin{array}{lll}
\delta H(t):= H(t,  x^*(t), u_i, m(t),{p}(t), {q}(t)) -H(t, x^*(t),   u^*_i(t), m(t),{p}(t) ,{q}(t)),
\end{array}
\end{equation}
\emph{and $H_k(t):=b_k(t)p+\sigma_k(t)q+r_k(t),$
for $k\in \{x,y, xx\}.$}
\end{proposition}

\section{Limitations and Challenges}
{\color{black}
The examples  above show that the continuum of agents assumption is rarely observed in engineering practice. The agents are not necessarily symmetric and a single agent may have a non-negligible effect on the mean field terms as illustrated in the HVAC application. 
 Without having a broad set of facts on which to theorize, there is a certain danger of mean-field game models that are mathematically elegant, yet have little connection to actual behavior observed in engineering practice. 
At present, our empirical knowledge is inadequate to the main assumptions of the classical mean-field game theory. This is why  a relaxed version is needed in order to better capture wide ranges of  behaviors and constraints observed in engineering systems. MFTG relaxations  include symmetry breaking, mixture between atomic and nonatomic agents, non-negligible effect on individual localized mean-field terms, and arbitrary  number of decision-makers. In addition,  behavioral and psychological factors should be incorporated for learning and information processes 
used by people-centric engineering systems. MFTG is still under development and  is far from being a well-established tool for engineered systems. 
Until now, MFTG was not focused on behavioral and cognitively-plausible models of choices in humans,  robots, machines, mobile devices and software-defined  strategic interactions. 
 Psychological and behavioral mean-field type game theories seem to explain behaviors that are better captured in experiments or in practice than classical game-theoretic equilibrium analysis. 
It allows  to consider psychological aspects of the agent in addition to the traditional "material" payoff modelling. The value  depends  upon choice consequences, mean-field states, mean-field actions and on beliefs about what will happen. 
The psychological MFTG framework can link cognition and emotion. It  expresses emotions, guilt, empathy, altruism, spitefulness (maliciousness) of the agents. 
It also include belief-dependent and other-regarding preferences in the motivations. It needs to be   investigated how much 
 the psychology of the people matters in  their behaviors in engineering MFTGs. The answer to this question is particularly crucial when analyzing the quality-of-experience of the users in terms of MOS (mean opinion score) values. A preliminary result  from a recent experiment conducted  in \cite{psy,innov} with 47 people carrying 
mobile devices with WiFi direct and D2D technology shows that the participation in forwarding the data of the users is correlated with their level of empathy towards their neighbors. This suggests the use of not only material payoffs but also non-material payoffs in order to better capture users behaviors. 
Another aspect of  MFTGs is the complexity of the analysis (both equilibrium and non-equilibrium) when multiple agents (and multiple mean-field terms) are involved in the interaction \cite{large,smc,cdc15,bou3,large2}.  
}
%
%\section{Introduction}
%
%Please use the AIMS template to prepare your manuscript, 
%before you submit to our journal. 
%Please read carefully the instructions for authors at
%http://www.aimspress.com \cite{authour1}. 
%These are important instructions and explanations. 
%Thank you for your cooperation.
%
%
%\section{Materials and Method}
%\subsection{Subheading}
%
%\subsubsection{Sub-subheading}
%The heading levels should not be more than 4 levels. 
%The font of heading and subheadings should be 12 point 
%normal Times New Roman. The first letter of headings 
%and subheadings should be capitalized.
%
%\section{Results}
%The body text is in 12 point normal Times New Roman, 
%the line space is at least 15 point.
% 
%\begin{table}[H]
%\begin{center}
%\caption{Caption of the table.}
%\begin{tabular}{ccc} \hline
% & & \\\hline
% & & \\
% & & \\
% & & \\\hline
% &(Table body should be created by MS word table function; three-line table is preferred.)
%\end{tabular}
%\end{center}
%\end{table}
%
%\begin{figure}[H]
%\begin{center}
%\includegraphics[scale=0.8]{figure.pdf}
%\caption{Legend of the figure.}
%\label{Fig1}
%\end{center}
%\end{figure}
%
%\begin{equation}
%  \text{[add an equation here; use MS Word or MathType equation function]}
%\end{equation}
%
%\section{Discussion}
%

\section{Conclusion and Future Work}
The article presented  basic applications of  mean-field-type game theory in  engineering, covering  key aspects such as  de-congestion in intelligent transportation networks, control  of virus spread over network, multi-level building evacuation,  next generation wireless networks, incentive-based demand satisfaction in smart energy systems, synchronization and coordination of nodes,  mobile crowdsourcing and cloud resource management. It appears from the wide ranges of applications and coverage that mean-field-type game theory is a promising  tool for engineering problems. However, the framework is still under development and needs to be  improved  to capture realistic behavior observed in practice. 
Possible extensions of the work described in this article include the study of  mean-field-type games  for risk engineering,  and 
an integrated mean-field-type game framework for smarter cities ranging from transportation to water distribution with ICT (Information Communication Technology), big data and human-in-the-loop among several other interesting directions.

\section*{Acknowledgments}This research work is  supported  by U.S. Air Force Office of Scientific Research under grant number FA9550-17- 1-0259.
The authors would like to thank the Editor and the anonymous reviewers for interesting and constructive comments on the manuscript.
The authors would like to thank the seminar participants at KTH Sweden for many extremely fruitful discussions and for their inputs on the first draft of the manuscript.

\section*{Conflict of Interest}
The authors declare no conflict of interest in this paper.

\end{document}